\documentclass[10pt,reqno]{amsart}

\numberwithin{equation}{section}

\usepackage{marginnote}

\usepackage[colorlinks=true,linkcolor=blue,citecolor=blue,linktocpage=true,urlcolor=black,hyperindex,breaklinks]{hyperref}
\usepackage{mathrsfs}

\usepackage{amssymb}
\usepackage{amsmath}
\usepackage{amsbsy}
\usepackage{amscd}
\usepackage{amsthm}
\usepackage{amsfonts}
\usepackage{enumerate}
    \usepackage[abbrev,nobysame]{amsrefs}
\renewcommand{\MR}[1]{} 

\usepackage{enumitem}

\usepackage{color}

\usepackage[capitalize,sort]{cleveref}
\newtheorem{theorem}{Theorem}[section]

\newtheorem{claim}[theorem]{Claim}
\newtheorem{lemma}[theorem]{Lemma}
\newtheorem{corollary}[theorem]{Corollary}

\newtheorem{proposition}[theorem]{Proposition}

\theoremstyle{definition}
\newtheorem{definition}[theorem]{Definition}

\newtheorem{remark}[theorem]{Remark}
\theoremstyle{definition}
\theoremstyle{remark}

\long\def\symbolfootnote[#1]#2{\begingroup\def\thefootnote{\fnsymbol{footnote}}
\footnote[#1]{#2}\endgroup}

\newcommand{\bfv}{\mathbf{v}}

\newcommand{\GL}{\operatorname{GL}}

\newcommand{\onto}{\xymatrix{\ar@{>>}[r]&}}
\newcommand{\da}[4]{\xymatrix{#1 \ar@<.5ex>[r]^{#2} \ar@<-.5ex>[r]_{#3} & #4}}

\newcounter{subconst}[subsection]

\newcounter{const}

\newcounter{CONST}

\newcommand{\id}{\mathrm{Id}}

\newcommand{\R}{\mathbb R}

\newcommand{\Z}{\mathbb {Z}}
\newcommand{\N}{\mathbb {N}}

\newcommand{\inv}{^{-1}}
\def\Lie{\mathrm{Lie}}

\def\Span{\mathrm{span}}

\newcommand{\Fol}{\mathcal{F}}
\def\fol{\Fol}
\newcommand{\td}{\tilde}
\newcommand{\wtd}{\widetilde}

  \newcommand{\restrict}[2]{{#1}{\restriction_{{ #2}}}}
  \def\sm{\smallsetminus}
  \def\calV{\mathcal V}
    \def\calE{\mathcal E}
      
    \def\calF{\mathcal F}
\def\calP{\mathcal P}
\def\calG{\mathcal G}
\def\calH{\mathcal H}
\title{Normal forms for contracting dynamics revisited}
\author[A.\ Brown]{
	Aaron Brown
}
\address{
	\parbox{0.5\textwidth}{
		Department of Mathematics\\
		Northwestern University\\}
}
\email{{awb@math.northwestern.edu}}

\author[A.\ Eskin]{
	Alex Eskin
}
\address{
	\parbox{0.5\textwidth}{
		Department of Mathematics\\
		University of Chicago\\
		5734 S University Ave\\
		Chicago, IL 60637\\}
}
\email{{eskin@math.uchicago.edu}}

\author[S.\ Filip]{
	Simion Filip
}
\address{
	\parbox{0.5\textwidth}{
		Department of Mathematics\\
		University of Chicago\\
		5734 S University Ave\\
		Chicago, IL 60637\\}
}
\email{{sfilip@math.uchicago.edu}}

\author[F.\ Rodriguez Hertz]{
	Federico Rodriguez Hertz
}
\address{
	\parbox{0.5\textwidth}{
		Department of Mathematics\\
		Pennsylvania State University\\
		328 McAllister Building\\
		State College, PA 16802\\}
}
\email{{fjr11@psu.edu}}

\begin{document}
\maketitle
\def\partP{\mathscr P}

\def\calA{\mathcal{A}}
\def\calB{\mathcal{B}}
\def\vol{\mathrm {vol}}
\def\e{\mathbf e}

 \begin{abstract}
We revisit the theory of normal forms for non-uniformly contracting dynamics.  We collect a number of lemmas and reformulations of the standard theory that will be used in other projects.  
\end{abstract}

\newcommand\iprod[1]{\langle#1\rangle}
\setcounter{tocdepth}{1}

\def\bfV{\mathbf V}

\def\bfA{\mathbf A}

\def\bfG{\mathbf G}
\def\bfU{\mathbf U}
\def\bfN{\mathbf N}
\def\bfL{\mathbf L}
\def\bfW{\mathbf W}

\renewcommand{\emph}[1]{{\bf{#1}}}

\def\1{\mathbf 1}
\def\supp{\mathrm{supp}}
\tableofcontents

\def\acts{\curvearrowright}



\section{Introduction}
The purpose of this note is to reformulate the theory of normal forms for non-uniform, non-stationary contracting dynamics.  This follows many such treatments in the literature including \cite{MR3642250,MR3893265,MR1908557,MR2090767}.
Our main interest is to the dynamics along (strong) stable manifolds for diffeomorphisms or coarse Lyapunov manifolds for actions of higher-rank abelian groups.  

We present  proofs of some auxiliary lemmas that already exist in the literature for completeness, especially those that are needed for other projects by the author.    For the main results concerning the existence, uniqueness, and measurability of the normal form coordinates in \cref{thm:normalforms}, we simply refer to the main result of  \cite{MR3642250}.  We formulate  a version of the linearization of dynamics  along stable manifolds (and more general invariant foliations with contracting dynamics) in \cref{sec7} below which is essentially a reformulation of the standard theory of normal forms.  A different treatment of this construction, with more details, will appear in a forthcoming project by the author together with Eskin, Filip, and Rodriguez Hertz.  

The outline of the paper is as follows:  In \cref{Sec1}, we present the main linear algebraic constructions to define subresonant polynomials.  In  \cref{sec2} we establish some growth estimates, we give a characterization of subresonant polynomials in terms of equivariance of sequences of suitable contracting dynamics.  In  \cref{sec3}, we repackage subresonant polynomial maps as linear maps in some extended vector space.  In  \cref{sec4} we explain the structures on manifolds.  In  \cref{sec5} we present lemmas about forwards regularity of cocycles.  In  \cref{sec6} we present a small lemma on temperedness.  Finally, in  \cref{sec7}, we assemble all the above to repackage the existence of normal form change of coordinates as a linear cocycle.  

\section{Subgraded linear algebra}
\def\SG{weight}
\def\SGs{weights}
\def\CSGs{Weights}

\renewcommand{\paragraph}{\subsection*}

\subsection{{\CSGs} on  vector spaces}\label{Sec1}  Let $V$   be a vector space over $\R$.  Unless otherwise stated all vector spaces are assumed to be finite dimensional.  
\subsubsection{\CSGs}
We have the following key definition.
\begin{definition} 
 \label{subgrad}  
 A \emph{{\SG}} 
 on  $V$ is a function 
$\varpi_V\colon V\to \R \cup \{-\infty\}$
such that for all $v,w\in V$, 
\begin{enumerate}
\item $\varpi_V(v) = -\infty$ if and only if $v=0$;
\item \label{SG222} $\varpi_V(c v) = \varpi_V(v)$ for all  $c\in \R\sm \{0\}$;
\item \label{SG333} $\varpi_V(v+w) \le  \max\{\varpi_V(v), \varpi_V(w)\}$.  
\end{enumerate}
\end{definition}
In most  applications, we will further  assume  $\varpi_V$ takes only negative values.  
Whenever clear from context, we  omit the  subscripts and write $\varpi$ for the implied {\SG}.

\subsubsection{Filtrations}
Let $V$ be equipped with a {\SG} $\varpi$.  Given $\lambda\in \R$ define subspaces $$V_{\le \lambda}:= \{ v\in V\mid \varpi(v) \le \lambda\}, \quad \quad V_{< \lambda}:= \{ v\in V\mid \varpi(v) < \lambda\}.$$
If  $V$ is  finite dimensional, then   $\varpi$ takes finitely many values $ \varpi(V\sm \{0\}) =  \{ \lambda_1, \dots, \lambda_\ell\}$, ordered so that $\lambda_1<\lambda_2 <\dots <\lambda_\ell$.  
  In this case, we also write  $V_i = \{ v\in V\mid \varpi(v) \le \lambda_i\}.$ 
We  also refer to $\lambda_1, \lambda_2, \dots, \lambda_\ell$ as the weights of $V$.  
The value $\dim (V_i/V_{i+1})$ is the \emph{multiplicity} of the weight $\lambda_i.$

At times we also write $\lambda_0 = -\infty$.

 We associate the  \emph{flag} $\calV_{\varpi_V}$ to a  {\SG} $\varpi_V$  on a finite-dimensional $V$ 
 consisting of nested subspaces: $$\calV_{\varpi_V}=\bigl\{\{0\}=V_{0}\subset V_1\subset\dots\subset V_\ell=V\bigr\}.$$
For each $1\le i\le \ell$, we obtain linear foliations $\calF_i$ of $V$ given by $\calF_i(v) = v + V_i.$
Note that $\calF_{i}(v)\subset \calF_{i+1}(v)$ for all $v\in V$.

\subsubsection{Adapted bases} Let $V$ be finite dimensional and equipped with a {\SG} $\varpi$ taking values $\lambda_1<\dots < \lambda_\ell$.  
\begin{definition}\label{def:adapbasis}An ordered basis $\{e_i\}$ of $V$ is \emph{adapted to $\varpi$} if for every $1\le i\le \ell$,
$$ V_{\le \lambda_i} = \Span\{e_1, \dots, e_{m_1 +\dots +m_i}\}.$$
A (unordered) basis  $\{e_i\}$ is  \emph{adapted to $\varpi$} if it has some order which makes it is adapted.  
\end{definition}
Note that \cref{def:adapbasis} implies  $\varpi(e_i)\le \varpi(e_j)$  for all $1\le i\le j\le \dim(V)$.   
In particular, this implies the following.  
\begin{claim}\label{claim:basiscomp}
Let $\{e_i\}$ be a basis of  $V$ \emph{adapted to $\varpi$}.  Given a non-zero $v\in V$ write $v= \sum c_i e_i$.  Then $\varpi(v)= \max \{\varpi (e_i) : c_i\neq 0\}.$
\end{claim}
\begin{corollary}\label{cor:easypoop}
Suppose $v = v_0 + \sum_{j=1} ^d v_j$ where $\varpi(v_j)<\varpi(v_0)$ for all $1\le j\le d$.  Then $\varpi(v) = \varpi(v_0)$.  
\end{corollary} 

\subsection{Induced {\SGs}}\label{sss:inducesSG}
Let $V$ and $W$ be finite-dimensional vector spaces equipped with {{\SGs}} $\varpi_V$ and $\varpi_W$, respectively.  We construct   {\SGs} on various associated vector spaces.  
\subsubsection{Subspaces} We equip a subspace $U\subset V$ with  the restriction of the  {\SG} $\varpi_V$ on $V$ to $U$.  
\subsubsection{Quotients}
Let $W$ be a subspace of $V.$  We equip the quotient $V/W$ with the {\SG} $$\varpi(v+W) = \min\{ \varpi _V(v+w): w\in W\}.$$
To see this defines a {\SG} we need only check the third axiom.  Let $v_1, v_2\in V$ and select $w_1,w_2\in W$ with $\varpi_V(v_1+W)=\varpi_V (v_1+w_1)$ and  $\varpi_V(v_2+W)=\varpi_V (v_2+w_2)$.
Then, 
\begin{align*}
\varpi(v_1+ v_2 + W)&= \min \{\varpi_V (v_1 + v_2 + w): w\in W\} \\
&\le \varpi_V(v_1 + v_2 + w_1 + w_2) \\&\le \max\{\varpi_V (v_1 + w_1), \varpi_V(v_2 +w_2)\}\\
&= \max\{\varpi (v_1 + W), \varpi(v_2 + W)\}.
\end{align*}

Suppose $V$ admits a direct sum decomposition $V= U\oplus W$ by subspaces $U,W\subset V$.  Equip  $U$ and $W$ with their restricted {\SGs}.   We say $U$ and $W$ are \emph{compatible} with $\varpi$ if for all $\lambda\in \R$, 
$$V_{\le \lambda}= U_{\le \lambda}\oplus W_{\le \lambda}.$$
\begin{lemma}\label{lem:compat}
Suppose $V= U\oplus W$ is a direct sum decomposition such that $U$ and $W$ are compatible with $\varpi$.  Then the identification of $V/W$ with $U$ preserves  {\SGs}.   That is, $\varpi(u+W)= \varpi(u)$ for all $u\in U.$
\end{lemma}
\begin{proof}
By definition,  $\varpi(u+W)\le \varpi(u+0)= \varpi(u)$.  On the other hand, set $\lambda= \varpi(u+W) $ and take  $w\in W$ with $ \varpi(u+w)=\lambda$.  
Since  $V_{\le \lambda}= U_{\le \lambda}\oplus W_{\le \lambda}$, there are $u'\in U_{\le \lambda}$ and $w'\in W_{\le \lambda}$ with $u+w=u'+w'$.  Since $V=U\oplus W$ is a direct sum decomposition, we have $u=u'$ and $w=w'$ whence $\varpi(u)\le \lambda$ and the conclusion follows.
%
%
\end{proof}

Note that the definition of compatibility is symmetric in $U$ and $W.$
The primary examples  of compatible subspaces that arise in the sequel are the following:    $U=V_{\le \lambda}$ and $W\subset V$ any subspace transverse to $U$ in $V$.  
Then, as vector spaces equipped with weights,  we naturally identify $V/W$ with  $U$  and $V/U$ with $W$.
 
\subsubsection{Tensor  and direct products}
On pure tensors in $V\otimes W$ define  $$\varpi(v\otimes w):= \varpi_V(v)+ \varpi_W(w).$$ 
This extends to a unique {{\SG}} on $V\otimes W$ using \cref{subgrad}\eqref{SG333}.  
For $k\ge 1$, we similarly obtain a {\SG}  $\varpi_{V^{\otimes k}}$ on $V^{\otimes k}$ given on pure tensors by $$\varpi_{V^{\otimes k}}(v_1\otimes\dots \otimes v_k) = \sum _{i=1}^k \varpi_V(v_i).$$
On the direct product $V\oplus W$, we  obtain a {\SG}   given by $$\varpi(v,w) = \max \{\varpi_V(v), \varpi_W(w)\}.$$

We make the following observations:   
Let $\{e_i\}$ and $\{f_j\}$ be bases adapted to $\varpi_V$ and $\varpi_W$ respectively.  Then 
\begin{enumerate}
	\item $\{e_i \otimes  f_j\}$ is a basis of $V\otimes W$ adapted to $\varpi_{V\otimes W}$, and 
	\item  $\{e_i\oplus 0\}\cup \{0 \oplus f_j\}$ is a basis of $V\oplus W$  adapted  $\varpi_{V\oplus W}$ 
\end{enumerate}

\subsubsection{Dual spaces}
Let $V^*$ denote the dual  of $V$.  Set $\varpi_{V^*} (0)= -\infty$.  Given a non-zero $\xi\in V^*$,  write   $$\varpi_{V^*}(\xi) = -\min\{ \varpi_V(v) : \xi (v) \neq 0\}= \max\{ -\varpi_V(v) : \xi (v) \neq 0\}.$$
One  verifies $\varpi_{V^*}$ defines  a {\SG} on $V^*$.  
If $\lambda_1<\lambda_2< \dots <\lambda_\ell$ denote the distinct values of $\varpi_V$ then $-\lambda_\ell < \dots <-\lambda_1$ are the distinct values of $\varpi_V^*$.  We similarly obtain a dual flag 
$$\{0\}= V^*_{-\infty} \subset V^*_{\le - \lambda_\ell}\subset \dots \subset V^*_{\le - \lambda_1} = V^*$$
where $V^*_{\le - \lambda_i} := \{\xi \in V^*: \varpi_{V^*}(\xi) \le - \lambda_i\}. $
We check, alternatively, the following characterization:
$$V^*_{\le - \lambda_i} = \{\xi \in V^*:V_{<\lambda_i} \subset  \ker (\xi)\}.$$
Indeed, if $v\in V_{<\lambda_i}$ and $\varpi(\xi ) \le -\lambda_i$
then $\xi(v)=0$ holds by definition.  Similarly, if $\xi(v) =0$ for all $v\in V_{<\lambda_i}$ then $\varpi(\xi)\le -\lambda_i$ holds by definition.


Equip $\R$ and $\R^*$ with the trivial {\SG} \begin{equation}\label{gradeonR}\varpi_\R(t) = \begin{cases}0 & t\neq 0\\ -\infty & t = 0.\end{cases}\end{equation}
With this definition we may alternatively  define the subgrading on $V^*$ by $$\varpi_{V^*}(\xi) = \max\bigl \{\varpi(\xi(v)) - \varpi(v): v\in V\sm\{0\}\bigr\}.$$
%
%
%
%

\subsection{Linear and multilinear maps}
Let $V$ and $W$ be finite-dimensional vector spaces equipped with {\SGs}  $\varpi_V$ and $\varpi_W$ taking values  $\lambda_1<\lambda_2 <\dots <\lambda_\ell$ and  $\eta_1<\eta_2 <\dots <\eta_p$, respectively. 
We  identify   $k$-multilinear maps $T\colon V^{\otimes k} \to W$ with   elements of $W\otimes (V^*)^{\otimes k}.$
The {\SGs} $\varpi_W$ and $\varpi_{V^*}$ on $W$ and $V^*$ then induce a {\SG}  $\varpi$ on the vector space of all $k$-multilinear functions $T\colon V^{\otimes k} \to W$.  

More explicitly, we have the following   characterization.
\begin{claim}\label{claim:maxmindef}
Let  $T\colon V^{\otimes k} \to W$ be a multilinear map.  Then    \begin{equation}\label{eq:wgtlinear}\varpi(T) = \max \left\{ \varpi_W\big(T(v_1, \dots ,v _k)\big) - \sum_{i=1}^k \varpi_V(v_k): v_1, \dots, v_k \neq 0\right\}.\end{equation}
\end{claim}

\begin{proof}
Let $\{e_j\}$ and $\{f_i\}$ be bases adapted to $\varpi_V$ and $\varpi_W$, respectively.  
Then $\{f_i  \otimes e_{j_1}^* \otimes \dots\otimes e_{j_k}^* : 1\le i\le p, 1\le j_\ell<k\}$ is a basis of $W\otimes (V^*)^{\otimes k}$ adapted to the {\SG} on $W\otimes (V^*)^{\otimes k}.$

We may assume $T \neq 0.$
Consider first the case $$T=  f_i \otimes e_{j_1}^* \otimes \dots\otimes e_{j_k}^*.$$
By definition,   $$\varpi( T )= \eta_i -\left( \lambda_{j_1} +\dots + \lambda_{j_k}\right).$$
Consider a non-zero pure tensor $\bfv= v_1\otimes  \dots\otimes  v_k\in V^{\otimes k}$.  We have 
$$\varpi_W\left (T(\bfv )\right) = \begin{cases} -\infty & \text{$   e_{j_\ell}^*(v_\ell)=0$ for some $1\le \ell \le k$}\\  \eta_i &\text{otherwise}. \end{cases}$$
In the second case that $e_{j_\ell}^*(v_\ell)\neq 0$ for all $1\le \ell \le k$, from \cref{claim:basiscomp} we have $\varpi_V(v_\ell)\ge \lambda_\ell$ for each $1\le \ell\le k$ whence
$$\varpi_{V^{\otimes k}} (\bfv) \ge  \lambda_{j_1} + \dots + \lambda_{j_k}$$ 
In particular, for all non-zero pure tensors $\bfv= v_1\otimes  \dots\otimes  v_k\in V^{\otimes k}$,
\begin{equation}\label{linsubadditive}\varpi_W(T(\bfv)) \le \eta_i\le \varpi(T) + \varpi_{V^{\otimes k}}(\bfv).\end{equation}
Moreover, with $\bfv = e_{j_1} \otimes \dots\otimes e_{j_k}$, we have 
$$\varpi_W(T(\bfv)) = \varpi(T) + \varpi_{V^{\otimes k}}(\bfv).$$
In particular,  \eqref{eq:wgtlinear} holds for $T=  f_i \otimes e_{j_1}^* \otimes \dots\otimes e_{j_k}^*$.

For general $T$, write $$T= \sum {}a_{i,j_1,\dots, j_k} f_i \otimes e_{j_1}^* \otimes \dots\otimes e_{j_k}^*.$$
Again by \cref{claim:basiscomp},
  \begin{align*}
  \varpi(T) &= \max \{\varpi(  f_i \otimes e_{j_1}^* \otimes \dots\otimes e_{j_k}^* ): a_{i,j_1,\dots, j_k} \neq 0 \}\\
 &= \max \left\{\eta_i -\left( \lambda_{j_1} +\dots + \lambda_{j_k}\right): a_{i,j_1,\dots, j_k} \neq 0  \right\}.
\end{align*}
Given a non-zero pure tensor $\bfv= v_1\otimes  \dots\otimes  v_k\in V^{\otimes k}$, from \eqref{linsubadditive} we have 
$$\varpi_W(T(\bfv)) -  \varpi_{V^{\otimes k}}(\bfv) \le \varpi(T) .$$
Let $(i,j_1,\dots, j_k)$ be so that $a_{i,j_1,\dots, j_k} \neq 0$
and $\varpi(T) = \eta_i -\left( \lambda_{j_1} +\dots + \lambda_{j_k}\right).$
Let $\bfv = e_{j_1} \otimes \dots\otimes e_{j_k}.$
Then  
$$ \varpi(T) = \varpi_{V^{\otimes k}}(\bfv)- \varpi_W(T(\bfv)) $$
and the claim follows.
\end{proof}
It follows that if $T\colon V\to W$ and $S\colon W\to U$ are linear maps then \begin{equation} \label{subadLin} \varpi (S\circ T) \le \varpi(S) + \varpi (T) .\end{equation}

We will say a linear map $T\colon V\to W$ is \emph{subresonant} if $\varpi(T)\le 0$ and \emph{strictly subresonant} if $\varpi(T)<0$.
In particular, a linear map $T\colon V\to W$ is {subresonant} if and only if $$\varpi_W(T(v))\le \varpi _V(v)$$ for all $v\in V$.

\subsubsection{Characterization and properties of invertible subresonant linear maps}
Suppose $V$ and $W$ are isomorphic  finite-dimensional vector spaces. 
\begin{definition}We say two {\SGs} $\varpi_V$ and $\varpi_W$ on $V$ and $W$ are \emph{compatible} if
\begin{enumerate}
\item the values taken by $\varpi_V$ and $\varpi_W$  coincide, and 
	\item the associated filtrations $\calV_{\varpi_V}$ and $\mathcal{W}_{\varpi_W}$ have the same multiplicities.
\end{enumerate}
\end{definition}
\begin{lemma}\label{lem:invert}
Let $V$ and $W$ have compatible {\SGs} and let $T\colon V\to W$  be a linear map with $\varpi(T)\le 0 $.  Then $T$ is invertible if and only if $\varpi_W(T(v)) = \varpi_V(v)$ for every $v\in V$.  
\end{lemma}

\subsection{Polynomial maps}
Given a degree $k\ge 1$ homogeneous polynomial $\phi\colon V\to W$, there exists a unique symmetric, $k$-multilinear function $L\colon V^{\otimes k} \to W$, called the \emph{polarization of $\phi$}, such that $$\phi(v) = L (v\otimes \dots \otimes v)$$ for all $v\in V$. 
When $k=1$, we have that $$\phi = L = D_0 \phi$$ coincides with the derivative of $\phi$ at $0$.  
If $k\ge 2$ then the polarization of $\phi$ is $$L= \frac{1}{k!} D^k_0\phi$$ where $D^k_0\phi$ is the total $k$th derivative tensor of $\phi$ at $0$.  
To avoid notational excess, given $k\ge 1$ and a $C^k$ function $f\colon V\to W$ we define $$ \wtd D^k f:=\frac {1} {k!} D^k_0 f.$$
Also let $\wtd D^0f= f(0)$.  

We remark that when  $k=0$,  we define a degree $0$ homogeneous function $\phi$ to be a constant function $\phi(v) = w_0$ for all $v\in V$; in this case, we have $V^{\otimes 0} = \R$ and we declare the polarization $L\colon \R\to W$ to be the linear map with $L(1) = w_0$ which we identify with  $$w_0 \otimes 1^*\in W\otimes \R^* = W\otimes (V^*)^{\otimes 0}.$$   
Given $w_1, \dots, w_k\in W$, under the canonical identification $W=W\otimes \R^*$, we similarly identify the following:
$w_1\otimes \dots \otimes w_k\otimes 1^*=
(w_1\otimes 1^*)\otimes \dots \otimes (w_k\otimes 1^*)
=w_1\otimes \dots \otimes w_k$.

\subsubsection{{\CSGs} and filtrations on polynomial maps}
For $k\ge 0$, let $S^k(V)$ denote the symmetric $k$-tensors on $V$; when $k=0$ we have $S^k(V) = \R$.   By slight abuse of notation, we identify a degree $d$ polynomial map $f\colon V\to W$ with an element of $$W \oplus \left(W\otimes S^1(V^*) \right)\oplus \dots\oplus \left(  W\otimes S^d(V^*) \right)\subset \bigoplus_{k=0}^\infty \left(W\otimes S^k(V^*)\right)$$
by identifying each degree $k\ge 0$ homogeneous component of $f$ with its polarization $\wtd D^k f$.  Write $$P(V,W)= \bigoplus_{k=0}^\infty \left(W\otimes S^k(V^*)\right)$$ for the (infinite-dimensional) vector space of all polynomial maps from $V$ to $W$.  
The {\SG} on each $W\otimes S^k(V^*)$ extends to a {\SG} $\varpi$ on the infinite-dimensional subspace of all polynomial maps $f\colon V\to W$.  
Under this identification, we equip the space $P(V,W)$ with the induced {\SG} $\varpi$ and write $$P_{\le \kappa}(V,W):= \{ f\colon V\to W: \varpi(f)\le \kappa\}$$   
for the induced filtration.
\begin{remark} 
Note that for a degree $k$ homogeneous polynomial  $f\in P(V,W)$ and $x\in V$, we have 
$$\varpi_W(f(x)) \le \varpi(f)  +  k \varpi_V(x).$$
Indeed, if $T=\wtd D^k f$ is the polarization of $f$, by \cref{claim:maxmindef},
$$\varpi_W(f(x))= \varpi_W(T(x,\dots, x)) \le   \varpi(T)  + \varpi_V(x^{\otimes k}) = \varpi(f)  + k \varpi_V(x).$$
However, for a degree $k\ge 2$ homogeneous function we do not in general have $$ \varpi(f)  = \max \{\varpi_W(f(x))- k \varpi_V(x): x\in V\}.$$
\end{remark}



We write $P(V) := P(V,\R) $ for the space of polynomial $\R$-valued functions on $V$.  Also write
$$P_{\le \kappa}(V) := P_{\le \kappa}(V,\R).$$
The vector space $P_{\le \kappa}(V)$ (for a suitably large $\kappa>0$) will be essential for later analysis.  
%

Given a polynomial  $f\colon V\to W$, write 
\begin{equation}
f= F_0 + F_1 + \dots + F_d\end{equation} in terms of the polarization of each homogeneous component where  $F_k = \wtd D^k f\in W\otimes S^k(V^*)$.  
If $V$ and $W$ are equipped with norms, let $\|F_k\|$ be the induced operator norm of each $F_k$ and write 
\begin{equation}\label{eq:normofpoly}\|f\|_P = \sum _{0\le k \le d} \|F_k\| =
 \sum _{0\le k \le d} \|
 \wtd D^k f\|\end{equation}
for the norm of the polynomial $f$.

\subsubsection{Polarizations under composition and subadditivitiy of \SGs}
Let $f\colon V\to W$ and $g\colon W\to U$ be polynomial functions.  
Write \begin{equation} f= F_0 + F_1 + \dots + F_d,  
\quad g= G_0 + G_1 + \dots + G_{d'}
\end{equation} where $F_i =\wtd D^if\colon V^{\otimes i}\to W$  and $G_j=\wtd D^jg \colon W^{\otimes j}\to U$ denote the polarization of the homogeneous components of $f$ and $g$, respectively.  
The degree $0$ component of $g\circ f$  has polarization 
\begin{equation}\label{polarizationzero} G_0 + \sum_{j = 1}^{d'} G_j\circ \left(F_{0}^{\otimes j}\right);\end{equation} 
for $k\ge 1$,  the degree  $k$ term of $g\circ f$ has polarization 
\begin{equation}\label{polarization}\sum_{j = 1}^{d'} \left(\sum_{i_1+ \dots +  i_j = k}G_j\circ (F_{i_1}\otimes  \dots \otimes F_{i_j})\right).\end{equation}
(If we declare the empty tensor product to be $1\otimes 1^*$, 
 we may similarly write \begin{equation}\label{polarization2}\sum_{j = 0}^{d'} \left(\sum_{i_1+ \dots +  i_j = 0}G_j\circ (F_{i_1}\otimes  \dots \otimes F_{i_j})\right).\end{equation}
 the polarization of the degree $0$  term of  $g\circ f$ in \eqref{polarizationzero}.)

From \eqref{polarizationzero}, \eqref{polarization}, and \eqref{subadLin} we obtain the following.  

\begin{proposition}[Subadditivity] \label{prop:subadd}  
Given polynomials functions $f\colon V\to W$ and $g\colon W\to U$ with $\varpi(f)\le 0$,
$$\varpi(g\circ f) \le \max\{\varpi(g(0)),  \varpi(g) + \varpi(f)\}.$$
\end{proposition}
\begin{proof}
We have $\varpi (F_{i_1} \otimes \dots \otimes F_{i_j}) \le \varpi (f)$ for all $1\le j\le 0$.  
\end{proof}

 \subsubsection{Families of subresonant polynomials} We collect a number of subsets of polynomials.     Let  $(V, \varpi_V)$ and $(W, \varpi_W)$ be vector spaces equipped with {\SGs}.
\begin{definition}
A polynomial map $f\colon V\to W$ is said to be
 \begin{enumerate} 
 \item  \emph{subresonant} if $\varpi(f) \le 0$;
\item \emph{weight decreasing} if $\varpi(f) \le  0$ and $\varpi(f(x))< \varpi(x)$ for all $0\neq x\in V$;
\item \emph{strictly subresonant} if $\varpi(f) < 0$.
 \end{enumerate}
We identify the following subvector spaces of $P(V,W)$: 
\begin{enumerate}
\item  $\calP^{SR}(V,W)$,  the space of   subresonant polynomial maps  $f\colon V\to W$; 
\item  $\calP^{SSR}(V,W)$, the   space of   strictly subresonant polynomial maps  $f\colon V\to W$;
\item  $\calP^{\ast}(V,W)$,  is the subspace of subresonant polynomial maps  $f\colon V\to W$ such that $\varpi(D_0f)<0$.
\end{enumerate}
\end{definition}

\begin{remark}\label{rem:SRvsWD} 
Write  $f\in P(V,W)$ as $f= F_0 + F_1 + \dots + F_d$ where $F_i$ is the polarization of the degree $i$ homogeneous component.  
  Then $\varpi(f) = \max \{\varpi (F_i)\}$; in particular, $f$ is subresonant (resp. strictly subresonant) if and only if each $F_i$ is subresonant (resp. strictly subresonant).

We clearly have $$\calP^{SSR}(V,W)\subset \calP^{*}(V,W)\subset \calP^{SR}(V,W).$$
However, we emphasize the first inclusion may be strict.  Indeed, consider $V=W=\R^2$ equipped with the {\SG} $\varpi(e_1) = -2, \varpi(e_2) = -1$.  Let $f\colon \R^2\to \R^2$ be $$f(x,y) = (y^2,0).$$
Then $\varpi(f)=0$ but $f\in \calP^{\ast}(V,W)$.  In particular, $f\in \calP^{\ast}(V,W)\sm \calP^{SSR}(V,W)$.   
\end{remark}

We observe the following: 
\begin{lemma}\label{weightdec}
Suppose $\varpi_V$ takes only negative values and let $f\in \calP^{SR}(V,W)$ be homogeneous of degree $k\ge 2$.  Then $f $ is weight decreasing.  
\end{lemma}
\begin{proof}
Let $T=\wtd D^k f$ be the polarization of $f$.  Since $\varpi(f)=\varpi(T)\le 0$ and $\varpi_V(x)<0$,  
$$\varpi_W(f(x)) = \varpi_W(T(x,\dots, x))\le \varpi_V (x,\dots, x) = k \varpi_V (x) <\varpi_V (x)$$
and thus $\varpi _W(f(x))< \varpi_V (x)$ for all $x\neq 0$.  
\end{proof}
As a corollary, any element of  $f\in \calP^{\ast}(V,W)$ (and thus $f\in \calP^{SR}(V,W)$) with $f(0) = 0$ is weight decreasing.

\subsubsection{Bounded degree and finite dimensionality}

\begin{lemma}\label{lem:bdddeg}
Suppose $\varpi_V$ takes only negative values $\lambda_1<\dots <\lambda_\ell<0$ and $\varpi_W$ take the values $\eta_1<\dots <\eta_p$.  Let $f\in P_{\le \kappa}( V, W)$.  Then $f$ has degree at most $\lfloor (\eta_1-\kappa)/\lambda_\ell \rfloor$.  In particular, $P_{\le \kappa}( V, W)$ is finite dimensional for all $\kappa \in \R$.   
\end{lemma}
\begin{proof}
Let $f\colon V\to W$ be homogeneous of degree $d$ with polarization $T\colon V^{\otimes d}\to W$.  If $f\in P_{\le \kappa}( V, W)$ then the image of $T$ is contained in $W_{\le (d \lambda_\ell +\kappa)}$.  If $d\lambda_\ell +\kappa\le \eta_1$ then $T \equiv 0$.  
\end{proof}

\subsubsection{Linearity and preservation of linear foliations}
\begin{lemma}\label{affineprops}
Suppose  $V$ and $W$ are equipped with compatible weights taking values $\lambda_1<\dots< \lambda_\ell<0$ and consider   $f\in \calP^{SR}(V,W)$.  
\begin{enumerate}
\item \label{affinish0} If $f$ is homogeneous of degree $k\ge 1$ then $f(V_{\le \lambda_i})\subset W_{\le \lambda_i}$.
\item \label{affinish1} For all $x\in V$, $f(x+V_{\le \lambda_i})\subset f(x) + W_{\le \lambda_i}$.  In particular, $f$   intertwines linear foliations of $V$ and $W$ parallel to the associated flags. 
\item \label{affinish3}  For each $1\le i\le \ell$, $f$ induces a subresonant polynomial map $\td f\colon  V/V_{\le \lambda_i}\to W/W_{\le \lambda_i}$;
\item \label{affinish2} If $f(0) = 0$ then the restriction  $f\colon V_{\le \lambda_1}\to W_{\le \lambda_1}$ is a linear map.
\end{enumerate}
\end{lemma}
\begin{proof}
\eqref{affinish0} 
follows from the definition of subresonance and \eqref{affinish1} 
follows immediately.
\eqref{affinish3}  follows from \eqref{affinish1}  and the definition of induced {\SG} on quotient spaces.  
\eqref{affinish2} follows from \cref{lem:bdddeg}
\end{proof}

\subsection{Properties of invertible subresonant polynomial maps}
We have the following standard fact about invertible subresonant polynomial maps.  We include a proof for completeness.  
\begin{proposition}\label{prop:invertible}
Suppose $V$ and $W$ have compatible {\SGs} taking only negative values.  Let $f\colon V\to W$ be a subresonant polynomial map such that  $D_0 f\colon V\to W$ 
is invertible.  Then 
\begin{enumerate}
	\item $f\colon V\to W$ is a diffeomorphism;
	\item the inverse function $f\inv \colon W\to V$ is a subresonant polynomial map.  
\end{enumerate}	
\end{proposition}
\begin{proof}
First consider the case that $f(0) = 0$.  
Write $f= F_1 + \dots + F_d$  where $F_i = \wtd D^if$ is the polarization of degree $i$ homogeneous component of $f$.  
We find a power series $g$ of the form 
 $$g = G_1 + G_2 + \dots$$ 
 where each $G_i\colon W^{\otimes i}\to V$ is symmetric.  
We argue $g$ is a formal inverse for $f$ and each $G_i$ is subresonant; it follows there are only finitely many $G_i$ and thus the power series converges.

Let $G_1= F_1\inv = (D_0 f)\inv.$
For $k\ge 1$, as in \eqref{polarization} the polarization of the degree $k$ term of $g\circ f$ is 
$$\wtd D^k(g\circ f)= \sum_{j = 1}^{k} \Bigl(\sum_{\substack{i_1+ \dots +  i_j = k\\  i_1, \dots, i_j \ge1}} G_j\circ (F_{i_1}\otimes  \dots \otimes F_{i_j})\Bigr) $$ and thus for $k\ge 2$ we aim to solve 
\begin{equation}\label{solvme}\sum_{j = 1}^{k} \Bigl(\sum_{\substack{i_1+ \dots +  i_j = k\\  i_1, \dots, i_j \ge1}} G_j\circ (F_{i_1}\otimes  \dots \otimes F_{i_j})\Bigr) =0.\end{equation}
Using \eqref{solvme} we may write $G_k (F_1 \otimes \dots \otimes F_1)$ in terms of $G_1, \dots, G_{k-1}$ and tensor powers of $F_j$.  
Since $F_1$ was assumed invertible, $F_1 \otimes \dots \otimes F_1\colon V^{\otimes k}\to V^{\otimes k}$ is invertible and we may inductively solve for $G_k$. By 
\cref{lem:invert}, $(F_1 \otimes \dots \otimes F_1)\inv$ is subresonant; moreover, by induction, each $G_k$ is thus subresonant.
Since the {\SG} is assumed to take only negative values, from \cref{lem:bdddeg}   only finitely many terms $G_k$ are non-zero; in particular, the formal power series for $g$ has only finitely many terms and so converges.  
We conclude that $g$ is a subresonant polynomial.  From \eqref{solvme}, $g\circ f (x) = x$ for all $x\in V$ and we conclude $g = f\inv$.  

Now suppose that $f(0) :=f_0\neq 0$.  If $g\colon W\to V$ is the inverse function of $x\mapsto f(x) -f_0$ then the map $$\hat g\colon x\mapsto g(x - f_0)$$ is a polynomial function with $\hat g\circ f (x) = x$ for all $x\in V$.  

From \cref{prop:subadd}$$\varpi(\hat g ) \le \max\{\varpi_W (g(0)), \varpi(g)+ 0\} \le 0 
$$
since  $g$ is subresonant.  It follows that $x\mapsto\hat  g(x)$ is subresonant.  
\end{proof}

\subsection{Properties under composition}
We have the following immediate corollary of \cref{prop:subadd}.   

\begin{corollary}
Fix $\kappa\in \R$, let $f\colon V\to W$ be a subresonant polynomial, and let $h\in P_{\le \kappa}(W,U)$.  Then $h\circ f\in  P_{\le \kappa}(V,U)$.   In particular, the map $h\mapsto  h\circ f$ induces a linear map between vector spaces
$$P_{\le \kappa}(W)\to P_{\le \kappa}(V).$$
\end{corollary}
\begin{proof}
Since $\varpi (f)\le 0$ and $\varpi (h)\le \kappa$, we have
$$\varpi(h\circ f) \le \max\{h(0), \varpi( h) + \varpi (f)\}\le \max \{\kappa, \kappa + \varpi(f)\} = \kappa. \qedhere$$
\end{proof}

We show that weights of non-constant terms are preserved by pre-composition by translation.

\begin{claim}\label{claim:translate}
Suppose $\varpi_V$ takes only negative values.  Write $f\in \calP^{SR}(V,W)$ as 	$f= F_0 + \dots + F_d$.  
Fix $v\in V$, let $\td f(x) = f(x + v)$,
and write 
$\td f= \td F_0 + \td F_1\dots + \td F_d$.  
Then for $1\le k\le d$,
 $$\varpi(F_k)= \varpi(\td F_k).$$
\end{claim}
\begin{proof}
For $k\ge 1$, the degree $k$ homogeneous component of $\td f$ is 
$$	\td F_k = \sum_{j=1}^d \sum_{i_1 + \dots + i_j = k} \left(F_j\circ ( Q_{i_1} \otimes \dots \otimes Q_{i_j})\right) $$
where $i_\ell\in \{0,1\}$ for  $1\le \ell\le j$,  $Q_0(x)=v$, and $Q_1(x)= x$.  
If $i_\ell =0$ for any $1\le \ell\le j$ then $j>k$, 
$$\varpi( Q_{i_1} \otimes \dots \otimes Q_{i_j}) <0,$$  and thus 
$$\varpi\left(F_j\circ ( Q_{i_1} \otimes \dots \otimes Q_{i_j})\right) < \varpi(F_j).$$ 
If $i_\ell=1$ for every $1\le \ell \le j$, then $j=k$ and 
$$F_j\circ ( Q_{i_1} \otimes \dots \otimes Q_{i_j}) = F_j.$$
It  follows from \cref{cor:easypoop} that $\varpi(\td F_k) =  \varpi (F_k).$
\end{proof}

\begin{lemma}\label{lemma:groupcomp}
Suppose $V$ is equipped with a {\SG} taking only negative values.  Let $P\subset \calP^{SR}(V,V)$ be either of  the vector subspaces: $P= \calP^{SSR}(V,V)$ or 
$P= \calP^{\ast}(V,V)$.  
Fix  $f,h\in \calP^{SR}(V,V)$ and $g\in P$. Then
\begin{enumerate}
\item \label{conc1} $g\circ f\in P$;
\item \label{conc2}$f\circ (h+g) = f\circ h + \td g$ where $\td g\in P$.
\end{enumerate}
\end{lemma}

\begin{proof}
Write \begin{equation} f= F_0 + F_1 + \dots + F_d,  \quad h= H_0 + H_1 + \dots + H_{d''}, \quad g= G_0 + G_1 + \dots + G_{d'}
\end{equation} 
in terms of polarizations of homogeneous components.  
We have $$\varpi(g\circ f) \le \max \{\varpi(g(0)), \varpi ( g) + \varpi (f)\}$$ and so conclusion  \eqref{conc1} holds if $P = \calP^{SSR}(V,V)$.  
Similarly, $D_0(g\circ f) = D_{f(0)}g \circ D_0f$.  By \cref{claim:translate}, if $g\in \calP^{\ast}(V,V)$ then $$\varpi (D_0(g\circ f) )\le \varpi (D_0(g)) + \varpi (D_0 f) <0.$$

For conclusion \eqref{conc2}, consider $k\ge 1$.  The degree $k$ homogeneous term of $f\circ (h+ g)$ has polarization 
$$\sum_{j = 1}^{d} \left(\sum_{i_1+ \dots +  i_j = k}\left(\sum_{
\delta_{i_1}, \dots \delta_{i_j}\in \{0,1\}}
F_j\circ (Q^{\delta_{i_1}}_{i_1}\otimes  \dots \otimes Q^{\delta_{i_1}}_{i_j})\right) \right)$$
where $Q^{0}_{i_j}= H_{i_j}$ and $Q^{1}_{i_j}= G_{i_j} $.
If $\delta_{i_\ell} = 0$ for all $1\le \ell \le j$,
then 
$$F_j\circ (Q^{\delta_{i_1}}_{i_1}\otimes  \dots \otimes Q^{\delta_{i_j}}_{i_j}) $$
contributes to 
$f\circ h$.
If $\delta_{i_\ell} = 1$ for some  $1\le \ell \le j$,
then  since each $G_j$ and $H_j$ are subresonant, 
$$\varpi(Q^{\delta_{i_1}}_{i_1}\otimes  \dots \otimes Q^{\delta_{i_j}}_{i_j} )\le \varpi(G_{i_\ell}).$$
and so $$\varpi\left(F_j\circ (Q^{\delta_{i_1}}_{i_1}\otimes  \dots \otimes Q^{\delta_{i_j}}_{i_j}\right)\le \varpi (G_{i_\ell}).$$
The conclusion of \eqref{conc2} then follows since if $g\in \calP^{SSR}(V,V)$ then $\varpi(G_\ell)<0$ for all $0\le \ell\le d'$ and 
if $g\in \calP^{\ast}(V,V)$ then $\varpi(G_\ell)<0$ for both $0\le \ell\le 1$ (which are the only terms of $g$ contributing to the degree 1 term of $f\circ (h+ g)$.)

\end{proof}

\subsection{Groups of invertible subresonant polynomial maps}
\def\Id{\mathrm{Id}}

Let   $V$ be    equipped with a {\SG} $\varpi_V$ taking negative values.  
We write  $$\calG^{SR}(V)\subset \calP^{SR}(V,V)$$ for the subset of all subresonant polynomials $f \colon V\to V$ such that $D_0f \colon V\to V$ is an invertible map.
We also write  $$\calG^{SSR}(V)\subset \calP^{SR}(V,V)$$ for the subset of all subresonant polynomials  of the form $$x\mapsto x + f (x)$$ where $f $ is strictly subresonant.  We refer to elements of $\calG^{SSR}(V)$ as \emph{invertible strictly subresonant maps.}
Write 
 $$\calG^{\ast}(V)\subset \calP^{SR}(V,V)$$ for the subset of all polynomial maps the subset of all subresonant polynomials  of the form $$x\mapsto x + f (x)$$ where $f\in \calP^{\ast}(V,V).$

\begin{claim} \ 
\begin{enumerate}
    \item \label{groups1}$\calG^{SR}(V)$ is a finite-dimensional Lie group.
    \item \label{groups2} Both $ \calG^{SSR}(V)$ and $ \calG^{\ast}(V)$ are normal subgroups of $ \calG^{SR}(V)$.  We have $ \calG^{SSR}(V)\subset  \calG^{\ast}(V)$.
    \item  \label{groups3} If $\varpi_V$ takes only negative values, $ \calG^{SSR}(V)$ contains all translations and thus acts transitively on $V$.
\end{enumerate}
\end{claim}
\begin{proof}
	Conclusion \eqref{groups1} follows from \cref{prop:subadd}   and \cref{prop:invertible}.
\cref{lemma:groupcomp}\eqref{conc1}  implies $ \calG^{SSR}(V)$ and  $\calG^{\ast}(V)$ are closed under multiplication.
\cref{lemma:groupcomp}\eqref{conc2} applied to $f\circ (f\inv + g\circ f\inv)$ implies $ \calG^{SSR}(V)$ and  $\calG^{\ast}(V)$ are fixed under conjugation by $\calG^{SR}(V)$.  It follows that  $ \calG^{SSR}(V)$ and  $\calG^{\ast}(V)$ are closed under inversion and are thus normal subgroups in $ \calG^{SR}(V)$.
\end{proof}

In \cref{uniprad}, we will see the following alternative characterization of $ \calG^{*}(V)$:
	$ \calG^{*}(V)$ is the unipotent radical of $ \calG^{SR}(V) $.

\section{Growth estimates and dynamical characterization of subresonant polynomial maps} \label{sec2}
\subsection{Growth estimates under composition}\label{sec:growth}
It will be useful to estimate the growth of  higher-order coefficients under repeated composition of subresonant polynomials.  When the  polynomials have linear part whose contraction rates are (nearly) as prescribed by the {\SG}, the coefficients of higher-order terms also decay as prescribed by the {\SG}.  

For each $j\in \N$, let  $V^j$ be a finite-dimensional vector space equipped with a {\SG}; we will assume all {\SGs} are all pairwise compatible and take values $$\lambda_1<\dots <\lambda_\ell<0.$$ Write $\varpi$ for the weight on each $V^j$.  
Fix $d = \lfloor \lambda_1/\lambda_\ell\rfloor\ge 1 $. 

 For each $j\ge 0$, 
 let 
  $$ f_j\colon V^j\to V^{j+1}$$
 be a  subresonant polynomial map
with  $  f_j(0)= 0$. 
   For $n\ge 1$, write $$f_0^{(n)} = f_{n-1} \circ \dots \circ f_0.$$

\begin{lemma}\label{lemgrowth}
Fix   $0<\epsilon $. 
Suppose there is $C\ge 1$ and a choice of norms on each $V^j$ 
such that the following hold:
\begin{enumerate}
\item \label{IPprop1} $ \|D_0f_jv\| \le e^{\varpi(v) +\epsilon}\|v\|$ for every $v\in V^j$,
\item \label{IPprop2}   $\|f_j\|_P \le C e^{\epsilon j}$.
\end{enumerate}

Then there exists $C_\epsilon\ge 1$ such  that for every $1\le k\le d$, every pure tensor ${\bfv} = v_1 \otimes \dots \otimes v_k \in (V^0)^{\otimes k}$, and every $n\ge 0$,
\begin{equation} \label{eq:growth2} 
\|\wtd D^kf^{(n)} _0\bfv\|\le C_\epsilon  e^{n(\varpi(\bfv) + 3k \epsilon)} 
\|\bfv\|. \end{equation}
\end{lemma}
\begin{proof}
We show the following stronger estimate: for $1\le k\le d$, there is $C_k\ge 1$
such  that 
for every pure tensor ${\bfv} = v_1 \otimes \dots \otimes v_k \in (V^0) ^{\otimes k}$ and $n\ge 0$,
\begin{equation} \label{eq:growth3} 
\|\wtd D^k  f_0^{(n)}\bfv\|\le (n+1)^{k-1}   C_k e^{n(\varpi(\bfv) + (2k-1) \epsilon)} 
\|\bfv\|. \end{equation}

Indeed in the case $k=1$, by hypothesis  we may take $C_1= 1$ and we proceed by induction on $k$.  Suppose for some $1\le k\le d-1$ and all $1\le i \le k$, we have found such $C_i$  satisfying \eqref{eq:growth3}.   
Set
\begin{align*}
	C_{k+1} =  C
	e^{-{(k+1) \lambda_1} } \sum_{j = 2}^{k+1}\left( \sum_{\substack{i_1+ \dots +  i_j = {k+1}\\  i_1, \dots, i_j \ge1}} 
  \left(
 \prod_{} C_{i_j} \right)\right) .
\end{align*}
When $n=0$, \eqref{eq:growth3} clearly holds since $C_{k+1} \ge 1$.   
We thus induct on $n$.  Assume \eqref{eq:growth3} holds for this $C_{k+1}$ and some $n$. 
We have 
\begin{align*}
\wtd D ^{k+1} f^{(n+1)} _0&= \sum_{j = 2}^{k+1} \left(\sum_{\substack{i_1+ \dots +  i_j = {k+1}\\  i_1, \dots, i_j \ge1}} \wtd D ^jf_n\circ \left(\wtd D ^{i_1}f^{(n)} _0\otimes  \dots \otimes \wtd D ^{i_j}f^{(n)} _0\right)\right) 
\\&\quad+  D_0 f_n \circ \wtd D ^{k+1}f^{(n)} _0
\end{align*}

Consider a pure tensor $\bfv \in (V^0)^{\otimes (k+1)}$.  We have \begin{equation}\label{bdddeg}\varpi (\bfv) \ge (k+1)\lambda_1. 
\end{equation}
 By the inductive hypothesis on $k$, 
\begin{align*}
\bigl \|&\wtd D ^{k+1} f^{(n+1)} _0\bfv\bigr\|\le \sum_{j = 2}^{k+1} \left(\sum_{\substack{i_1+ \dots +  i_j = {k+1}\\  i_1, \dots, i_j \ge1}} \bigl\|\wtd D ^jf_n\circ \bigl(\wtd D ^{i_1}f^{(n)} _0\otimes  \dots \otimes \wtd D ^{i_j}f^{(n)} _0\right)\bfv  \Bigr\| 
\Bigr) 
\\&\quad +\Bigl\| \wtd D ^{1}f_n \circ \wtd D ^{k+1}f^{(n)} _0\bfv \Bigr\|\\
&\le  \sum_{j = 2}^{k+1}\left( \sum_{\substack{i_1+ \dots +  i_j = {k+1}\\  i_1, \dots, i_j \ge1}} 
C e^{n \epsilon }
\Bigl( \prod_{} C_{i_j} \Bigr) (n+1)^{\sum(i_j-1)}  e^{n(\varpi(\bfv) +  {\sum(2i_j-1)} \epsilon)}  \|\bfv\|
  \right)  \\
&\quad   + e^{\varpi (\wtd D ^{k+1}f^{(n)} _0\bfv) + \epsilon } \Bigl\| \wtd D ^{k+1}f^{(n)} _0\bfv \Bigr\|\\
&\le  \sum_{j = 2}^{k+1} \left( \sum_{\substack{i_1+ \dots +  i_j = {k+1}\\  i_1, \dots, i_j \ge1}} 
C e^{n \epsilon }
\Bigl( \prod_{} C_{i_j} \Bigr)
 (n+1)^{k+1-j}  e^{n(\varpi(\bfv) +  {(2k+ 2 -j)} \epsilon)}  \|\bfv\|
  \right)\\
&\quad   + e^{\varpi (\bfv) + \epsilon } \Bigl\| \wtd D ^{k+1}f^{(n)} _0\bfv \Bigr\|\\
&\le  \Bigl( C e^{ {-(k+1)\lambda_1}}\sum_{j = 2}^{k+1} \sum_{\substack{i_1+ \dots +  i_j = {k+1}\\  i_1, \dots, i_j \ge1}} 
\Bigl( \prod_{} C_{i_j} \Bigr)
 \Bigr) (n+1)^{k-1}  e^{(n+1)(\varpi(\bfv) +  {(2k+1)} \epsilon)}  \|\bfv\|
\\
&\quad   + e^{\varpi (\bfv) + \epsilon } \bigl\| \wtd D ^{k+1}f^{(n)} _0\bfv \bigr\| \\ 
&= C_{k+1} (n+1)^{k-1}  e^{(n+1)(\varpi(\bfv) +  {(2k+1)} \epsilon)}  \|\bfv\|
  + e^{\varpi (\bfv) + \epsilon } \bigl\| \wtd D ^{k+1}f^{(n)} _0\bfv \bigr\|.
\end{align*}
Recall the inductive hypothesis  on $n$:
 $$\bigl\| \wtd D ^{k+1}f^{(n)} _0\bfv \bigr\| \le C_{k+1} (n+1)^{k}  e^{n(\varpi(\bfv) +   {(2k+1)} \epsilon)}  \|\bfv\|.$$
Then 
 \begin{align*}
\|\wtd D ^{k+1} &f^{(n+1)} _0\bfv\|\\
&\le C_{k+1}   e^{(n+1)(\varpi(\bfv) +   {(2k+1)} \epsilon)}  \Bigl((n+1)^{k-1} +(n+1)^{k}   \Bigr)\|\bfv\|\\
&\le C_{k+1}   e^{(n+1)(\varpi(\bfv) +   {(2k+1)} \epsilon)}  (n+2)^{k}   \|\bfv\|\\
&= C_{k+1}   e^{(n+1)(\varpi(\bfv) +   {(2(k+1)-1)} \epsilon)} ((n+1)+1)^{(k+1)-1}    \|\bfv\|.
\end{align*}
Claim \eqref{eq:growth3} then follows for all $1\le k\le d$, all $n\ge 0$, and all pure tensors $\bfv$.  

Finally, we take $$C_\epsilon =\sup\Bigl \{e^{-\epsilon n} C_k (n+1)^{k-1} :1\le k\le d, n\ge 0\Bigr\}$$
and conclude \eqref{eq:growth2}.  
\end{proof}


Applying \cref{lemgrowth} to pure tensors of the form $x\otimes \dots \otimes x$, we immediately obtain the following.
\begin{corollary}\label{corgrowth}
Fix   $0<\epsilon. $ 
Suppose the sequence of norms on $V^j$ are as in \cref{lemgrowth}.  
Then  for all $n$ and $x\in V^0$,   
 \begin{equation} \label{eq:growth2} 
\bigl\|D_0f^{(n)}_0(x)- f^{(n)} _0(x)\bigr\|\le d C_\epsilon  e^{n (2 \varpi (x) + 3d \epsilon)}\max\{ \|x\|^2, \|x\|^{d}\}\end{equation}
\end{corollary}
This has the following consequence.
\begin{corollary}\label{corgrowth2} ~
Assume $0<\epsilon< \frac{-\lambda_\ell }{10d } $.
\begin{enumerate}[label=(\alph*), ref=(\alph*)]
\item \label{parta}
If the norms on $V^j$ satisfy \eqref{IPprop1} and \eqref{IPprop2} of  \cref{lemgrowth} then for every compact $K\subset V_0$ there is $C_K\ge 1$ such  that 
for every $x\in K$  and   $n\ge 0$,
\begin{equation}\label{growthupper}\|f_0^{(n)}(x)\| \le C_K  e^{n(\varpi_0(x) + \epsilon)}\|x\|.\end{equation}

\item Suppose  in addition to \eqref{IPprop1} and \eqref{IPprop2} of \cref{lemgrowth}, the norms on $V^j$ satisfy the following:
for every $v\in V^0$ there is $n_0=n_0(v)$ so that for all $n\ge n_0$, $$ \|D_0f_0^{(n)}v\| \ge e^{n(\varpi_0(v)-\epsilon)}\|v\|$$
Then for every $x\in V^0$ there is $C_x\ge 1$ so that for all $n$, 
\begin{equation} \label{growthlower}\|f_0^{(n)} (x)\| \ge \frac 1 {C_x} e^{n(\varpi_0(x) -\epsilon)}\|x\|.\end{equation} 
\end{enumerate}
\end{corollary}

\subsection{Smooth equivariant functions are subresonant polynomials}
As in  \cref{sec:growth}, for each $j\in \N$, let  $V^j$ and $W^j$ be a finite-dimensional normed vector spaces equipped with a {\SG}.   We assume all  $V^j$  (resp.\ all $W^j$) are isomorphic and the {\SGs} are all pairwise compatible.  
Write $\varpi$ for the {\SG} on all spaces taking values $\lambda_1<\dots< \lambda_\ell <0$ on $V^j$ and $\eta_1<\dots< \eta_p <0 $ on $W^j$.
Let $$d_1:=\lfloor \lambda_1/\lambda_\ell \rfloor, \quad \quad d_2:=\lfloor \lambda_1/\eta_p \rfloor, \quad \quad d_3:=\lfloor \eta_1/\eta_p \rfloor,$$
and $d= \max\{d_1, d_2, d_3\}.$

Write $W^j(\rho)$ for the ball of radius $\rho$ centered at $0$ in $W^j$.  Given a $C^r$ function $\phi\colon W^j(\rho)\to V_j$,   write $\|\phi\|_{C^r}$ for the usual $C^r$ norm.  

Set  $$\epsilon_0 = \frac{1}{10d} \min \{1,-\lambda_\ell,  -\eta_p, \lambda_i-\eta_{j}:  \lambda_i\neq\eta_{j} \}.$$

\begin{proposition}\label{prop:forcedSR}
For every $j$, let $ f_j\colon V^j\to V^{j+1}$ and $ g_j\colon W^j\to W^{j+1}$ be  invertible subresonant polynomials with $f_j(0) = g_j(0) = 0$. 
 Fix any $r>\lambda_1/\eta_p$ and $0<\epsilon <\min \left\{\epsilon_0, \frac {\lambda_1 - d_2\eta_p}{3d_2+2}\right\}$.  
Suppose $V^j$ and $W^j$ are equipped with norms satisfying the  following:
\begin{enumerate}
\item  $ \|D_0f_jv\| \le e^{\varpi(v) +\epsilon}\|v\|$ and $ \|D_0g_jw\| \le e^{\varpi(w) +\epsilon}\|w\|$ for every $v\in V^j$ and $w\in W^j$,
\item   $\|f_j\|_P, \|g_j\|_P \le C e^{\epsilon j}$
\item \label{prop3norms} for every  $v\in V^0$, there is $n_0(v)$ such that for all $n\ge n_0(v)$,  $$e^{n(\varpi(v)-\epsilon )} \|v\| \le \|D_0f_0^{(n)}v\|. $$
\end{enumerate}

 For every  $n$, let $$\phi_n\colon W^n(\rho) \to V^n$$ be a $C^r$ function. 
Suppose the following hold:
\begin{enumerate}[resume]
\item $f_0^{(n)}\circ \phi_0(x) = \phi_n \circ g_0^{(n)}(x)$ for every $n$  and  every $x\in W^0(\rho)$;
\item \label{hyp5} $\sup _{j} \|\phi_{n_j}\|_{C^r}<\infty$  for some infinite subset $\{n_j\}\subset \N$.
\end{enumerate}

Then $\phi_0\colon W^0(\rho)\to V^0$ coincides with the restriction to $W^0(\rho)$ of  a subresonant polynomial map $h\colon W^0\to V^0.$
\end{proposition}

\begin{proof} 
For each $n$, write $$\phi_n = h_n + R_n$$ where $h_n\colon W^n\to V^n$ is a degree $d_2$ 
polynomial.  
Note there is $C(d_2)\ge 1$ so that  $$\|h_n\|_P\le C(d_2) \|\phi_n\|_{C^r}.$$  Fix $s$ with $d_2<s\le \min \{r, d_2+1\}$.   Then  
  $R_n\colon W^n\to V^n$ satisfies  $$\|R_n (x) \|\le \|\phi_n\|_{C^r} \|x\|^s$$ for all $x\in W^0(\rho)$.
For each $0\le k\le d_2$,  uniqueness of $k$-jets implies \begin{equation}\label{eqtaylor}\wtd D ^k(h_n \circ g_0^{(n)} )= \wtd D ^k(f_0^{(n)}\circ h_0)\end{equation}
and \begin{equation}\label{eqtaylor2}\wtd D ^k(h_n) = \wtd D ^k\bigl(f_0^{(n)}\circ h_0\circ( g_0^{(n)} )\inv\bigr).\end{equation}

\subsection*{Step 1: Each $h_n $ is subresonant}We  first claim that   $h_0\colon W^0\to V^0$ is subresonant polynomial map.  Note then that \eqref{eqtaylor2}  implies  each $h_n$ is subresonant and, using that both sides of \eqref{eqtaylor2} are subresonant and thus have degree at most $d_2$,    for all $n$ we have \begin{equation}\label{eqtaylor3}h_n \circ g_0^{(n)} = f_0^{(n)}\circ h_0.\end{equation}  

Suppose $h_0\colon  W^0\to V^0$ is not subresonant.   We have $\wtd D^0 h_0 = h_0(0) \in V^0$ and so $\varpi(\wtd D^0 h_0)\le \lambda_\ell<0$ and $\wtd D^0 h_0$ is subresonant.  
Let $1\le k \le d_2$ be the minimal degree for which $\wtd D^k  h_0$ is not  a subresonant multilinear map.  
We have 
$$\wtd D^k  (f_0^{(n)} \circ h_0 )= \sum_{j = 1}^{d_1} \left(\sum_{\substack{i_1+ \dots +  i_j = k\\  i_1, \dots, i_j \ge 0}} \wtd D^j  f_0^{(n)}\circ (\wtd D ^{i_1}h_0\otimes  \dots \otimes \wtd D ^{i_j}h_0)\right).$$
Fix a pure tensor $\bfv\in (W^0)^{\otimes k}$ for which $$\varpi(\wtd D^k  h_0 (\bfv))> \varpi (\bfv).$$
Set $\eta = \varpi (\bfv)$ and $\kappa = \varpi(\wtd D^k  h_0 (\bfv))$.
By the choice of $\epsilon <\epsilon_0$, we have $$\eta + 3d_1 \epsilon \le \kappa - 7 \epsilon, \quad \quad \kappa +\lambda_\ell + 3d_1 \epsilon \le \kappa - 7\epsilon.$$

We have 
\begin{align*}
	\wtd D^k(f_0^{(n)} \circ h_0 )(\bfv)
	&= 
	\sum_{j = 2}^{d_1} \left(\sum_{\substack{i_1+ \dots +  i_j = k\\  i_1, \dots, i_j \ge 0}} \wtd D^jf_0^{(n)}\circ (\wtd D^{i_1}h_0\otimes  \dots \otimes \wtd D^{i_j}h_0)\right)(\bfv)\\
	& \quad + D_0f_0^{(n)} \circ \wtd D ^kh_0 (\bfv).
\end{align*}
Recall that    $\wtd D^{i}h_0$ is subresonant for $1\le i \le k-1$ and that
 and $\wtd D^0 h_0 := h_0(0) \in V^0 $ and so $\varpi(\wtd D^0 h_0)\le \lambda_\ell<0$.  For any $j\ge 2$ and $i_1, \dots, i_j \ge 0$ with $i_1+ \dots +  i_j = k$, we have 
\begin {align*}
\varpi \bigl (&\wtd  D^{i_1}h_0\otimes  \dots  \otimes \wtd D^{i_j}h_0)\bigr)(\bfv) 
\\ &\le \begin{cases}
\eta & i_1, \dots, i_j \ge1\\
  \eta + \lambda_\ell & \text{$i_p= 0$ for some $1\le p\le j$ and $i_q<k$ for all $1\le q\le j$}\\
  \kappa + \lambda_\ell & \text{$i_p= 0$ for some $1\le p\le j$ and $i_q=k$ for some $1\le q\le j$}.
\end{cases}
\end{align*}
In all cases, for $j\ge 2$ we have $$\varpi \bigl (\wtd  D^{i_1}h_0\otimes  \dots  \otimes \wtd D^{i_j}h_0)\bigr)(\bfv) \le  \kappa + \lambda_\ell$$
and so 
\begin{equation}\label{eq:feces}
\varpi \left (\wtd D^{i_1}h_0\otimes  \dots \otimes \wtd D^{i_j}h_0)\right)(\bfv) +3d_1 \epsilon \le \kappa - 7 \epsilon. 
\end{equation}

By \cref{lemgrowth}  and \eqref{eq:feces}, there is $C_1\ge 1 $ such  that for all $n\ge 0$ we have the upper bound 
\begin{equation}
\left\|		\sum_{j = 2}^{d_1} \left(\sum_{\substack{i_1+ \dots +  i_j = k\\  i_1, \dots, i_j \ge 0}} \wtd D^jf_0^{(n)}\circ (\wtd D^{i_1}h_0\otimes  \dots \otimes \wtd D^{i_1}h_0)\right)(\bfv)\right\|
\le  C_1 e^{n(\kappa - 7 \epsilon)} 
\|\bfv\|. \label{flow}
\end{equation}
Since $\wtd D ^kh_0 (\bfv)\neq 0$, by property \ref{prop3norms} of the norms on $V^j$, for all sufficiently large $n$ we have a lower bound 
\begin{equation}
	\left\| D_0f_0^{(n)} \circ \wtd D ^kh_0 (\bfv)\right\|
	\ge  e^{n(\kappa -2\epsilon)} \label{fupper}
\end{equation}
On the other hand, \cref{corgrowth2}\ref{parta} implies there is a $C_2$ such that for all $n\ge 0$ we have the upper bound
\begin{align}
\|\wtd D^k  (h_n\circ g_0^{(n)} ) \bfv\| & = \left \| \sum_{j = 1}^{k}\left(\sum_{\substack{i_1+ \dots +  i_j = k\\  i_1, \dots, i_j \ge1}} \wtd D^j h_n\circ (\wtd D^{i_j} _0g_0^{(n)}\otimes  \dots \otimes \wtd D^{i_1}_0g_0)\right)(\bfv)	\right \| \notag \\
&\le C_2 ( \|h_n\|_P) e^{n(\eta + 3d_2 \epsilon)} \|\bfv\|\\\
&\le   C(d_2) C_2  \|\phi_n\|_{C^k} e^{n(\kappa - 7 \epsilon)} \|\bfv\|\label{glow}
\end{align}
Since we assume $ \|\phi_{n_j}\|_{C^k} $ is uniformly bounded for some infinite set $\{n_j\}$, the bounds in \eqref{flow}, \eqref{fupper}, and \eqref{glow} contradict the equality of degree $k$ terms in  \eqref{eqtaylor}.

\subsection*{Step 2: $R_0\equiv 0$.}  
Having shown each $h_j$ is subresonant, 
we conclude  $R_0\equiv 0$.  

For the sake of contradiction, suppose there is $x\in W^0(1)$ with $$\phi_0(x) -h_0(x) = R_0(x) \neq 0.$$ 
Let $y_0 = h_0(x)$ and $y_1 = R_0(x)$.
We have
$$f_0^{(n)} \bigl(\phi_0(x) \bigr)= f_0^{(n)} \bigl(y_0 +y_1 \bigr)= \sum _{j=1}^d \left(\sum_ {\substack{\delta_{i}\in \{0,1\}}}
\wtd D^jf _0^{(n)}(y_{\delta_1}, y_{\delta_2},\dots, y_{\delta_j})\right).
$$
Observing that  terms with $\delta_i=0$ for all $1\le i\le j$ contribute to $f_0^{(n)}  \bigl( h_0(x) \bigr) $ 
and terms with 
 $\delta_i=1$ for all $1\le i\le j$ contribute to $f_0^{(n)}  \bigl( R_0(x) \bigr) $, we write  
$$f_0^{(n)}(\phi_0(x)) = f_0^{(n)}  \bigl( h_0(x) \bigr) + f_0^{(n)}  \bigl( R_0(x) \bigr)  
+ \td r_n(x)$$
where 
$$\td r_n(x)=\sum _{j=2}^d \left(\sum_ {\substack{\delta_{i}\in \{0,1\}\\ 1\le \sum \delta_i\le j-1}}
\wtd D^jf_0^{(n)}(y_{\delta_1}, y_{\delta_2},\dots, y_{\delta_j})\right).
$$
By \cref{lemgrowth}, there is $C$ independent of $n$ and $x$ such that
\begin{equation}\label{eq:errorgrowth}
\|\td r_n(x)\| \le C  e^{n(\varpi(h_0(x) ) + \varpi(R_0(x)) + 3d \epsilon)}.	
\end{equation}

We have 
\begin{align*}
h_n \circ   \bigl( g_0^{(n)}(x) \bigr)  &+ R_n \circ   \bigl( g_0^{(n)}(x) \bigr)\\
&= \phi_n\circ  \bigl( g_0^{(n)}(x) \bigr)\\
&= 
f_0^{(n)} \circ  \phi_0(x)\\
&=f_0^{(n)}  \bigl( h_0(x) + R_0(x)\bigr) \\
&=f_0^{(n)}  \bigl( h_0(x) \bigr) + f_0^{(n)}  \bigl( R_0(x) \bigr) + \td r_n(x).
\end{align*}
Since $h_n \circ   \bigl( g_0^{(n)}(x) \bigr)= f_0^{(n)}  \bigl( h_0(x) \bigr) $, we   have
 \begin{equation}\label{eq:dungeater}f_0^{(n)} \bigl(R_0(x)\bigr)  = R_n \circ   \bigl( g_0^{(n)}(x) \bigr)-  \td r_n(x).\end{equation}

By  \cref{corgrowth2}, for all $n\ge 0$ sufficiently large,
\begin{equation}\label{eq:lowbd}
	\| f_0^{(n)} \bigl(R_0(x)\bigr) \| \ge e^{n(\varpi(R_0(x)) -2\epsilon)} 
\end{equation}
and
\begin{equation} \label{eq:upbd}
	\left\|R_n\left(g^{(n)}_0 (x)\right)\right\|\le \|\phi_n\|_{C^r} e^{ns(\eta_p+ 2\epsilon) }.
\end{equation}
By the choice of $\epsilon>0$, we have 
$$\varpi(R_0(x)) -2\epsilon > \lambda_1-2\epsilon >s\eta_p+ 3s\epsilon$$
and
$$\varpi(R_0(x)) -2\epsilon >  \varpi(h_0(x) ) + \varpi(R_0(x)) + 3d  \epsilon.$$
Using that $\|\phi_{n_i}\|C^{r}$ is bounded for  some infinite subset of $n_i$, the estimates \eqref{eq:errorgrowth},  \eqref{eq:lowbd}, and \eqref{eq:upbd} contradict equality in \eqref{eq:dungeater}.
\end{proof}

\section{Linearization of subresonant polynomials}\label{sec3}
Let $V$ and $W$ be finite-dimensional vector spaces equipped with compatible {\SGs} taking negative values. By passing to  suitable quotients of tensor powers of $V$ and $W$,  every subresonant polynomial map $f\colon V\to W$ is canonically  identified with a linear map between associated vector spaces.

\subsection{Linearization theorem}
\begin{theorem}[Linearization of subresonant polynomials]\label{thm:LinSRPoly}
Let $V$ be a finite-dimensional vector space equipped with a  {\SG} $\varpi_V$ taking only negative values.

There exists a finite-dimensional vector space $\bfV$ equipped with a weight $\varpi_{\bfV}$ taking non-positive values, a linear map $\Pi_V\colon \bfV\to V$, and a polynomial map $\iota _V\colon V\to \bfV$ 
 with the following properties:
\begin{enumerate}[label=(\arabic*),ref=(\arabic*)]
\item \label{LineNF0}
$\bfV':=\{v\in \bfV: \varpi_{\bfV}(v)<0\}$ is a a codimension-1 subspace of $\bfV$. 
\item \label{LineNF1} The image of $\iota_V$ is contained in a  codimension-1 affine subspace   $\iota_V(0)+ \bfV'$. 
\item \label{LineNF2} For every $f \in \calP^{SR}(V)$ there exists a   linear map $\bfL f \colon \bfV\to \bfV$ such that $$\iota_V\circ f    = (\bfL f )\circ \iota_V.$$
\item \label{LineNF3} Let $\bfG^{SR}(\bfV)$,  $\bfG^{SSR}(\bfV)$,  and $\bfG^{\ast}(\bfV)$ denote the images of  $\calG^{SR}(V)$,  $\calG^{SSR}(V)$, and  $\calG^{\ast}(V)$, respectively,  under  $f \mapsto \bfL f $.  Then   $ \bfG^{SR}(\bfV)$, $\bfG^{SSR}(\bfV)$, and $\bfG^{\ast}(\bfV)$ are subgroups of $\GL(\bfV)$ and the map $f \mapsto \bfL f $ determines a continuous isomorphism between Lie groups.  
\item \label{LineNF4}  $\bfG^{\ast}(\bfV)$ and $\bfG^{SSR}(\bfV)$ are unipotent subgroups of $\bfG^{SR}(\bfV)$.  
\item \label{LineNF5} $\Pi_V\circ \iota_V$ is the identity map.  
\end{enumerate}
Let $U$ and $W$ be finite-dimensional vector spaces equipped with   {\SGs} compatible with the {\SG} on $V$.  
\begin{enumerate}[label=(\arabic*),ref=(\arabic*),resume]
\item \label{LineNF7} For every $f  \in \calP^{SR}(V,W)$ there is a   linear transformation $\bfL f \colon \bfV\to \bfW$ with $$\iota_W\circ f  = (\bfL f )\circ \iota_V.$$
It follows  $\bfL f(\bfV')\subset \bfW'$.
\item \label{LineNF8} Given  $f  \in \calP^{SR}(V,W)$ and $g \in \calP^{SR}(W,U)$,
	$$\bfL(g\circ f ) = (\bfL g)\circ (\bfL f ).$$
\end{enumerate}
Let $\bfv_0= \iota_V(0)$ and let $\ell_V$ denote the span of $\{\bfv_0\}$.  
\begin{enumerate}[label=(\arabic*),ref=(\arabic*),resume]
\item \label{LineNF6} An inner product on $V$ canonically defines an inner product on $\bfV$; relative to this inner product, $\iota_v(0)$ is a unit vector, $\ell_V$ is orthogonal to $\bfV'$, and $D_0\iota_V$ is an isometry onto its image.    
\item \label{LineNF9} If $f  \in \calP^{SR}(V,W)$ satisfies $f (0) = 0$ then $\bfL f \colon \bfV\to \bfW$ induces a linear map between the quotients
$$\bfV/\ell_V\to \bfW/\ell_W.$$
Using the orthogonal projections $\bfV\to \bfV'$ and  $\bfW\to \bfW'$ along $\ell_V$ and  $\ell_W$, we may view this as a map $\bfV'\to \bfW'$.  
\end{enumerate}
\end{theorem}

\subsection{Construction of the linearization}  
Let $V$ and $W$ have compatible {\SGs} taking values $\lambda_1<\dots< \lambda_\ell<0.$
Given any $f \in \calP^{SR}(V,W)$  and $\kappa\in \R$, by \cref{prop:subadd} we obtain a linear map $f ^*\colon P_{\le \kappa}(W) \to P_{\le \kappa} (V)$ given by precomposition: $$f ^*(h) =h\circ f .$$

Fix once and for all  $\kappa= -\lambda_1$.
We have that elements of $P_{\le - \lambda_1}(V)$ and $P_{\le -\lambda_1}(W)$ have degree at most $d_0=\lfloor \lambda_1/\lambda_\ell\rfloor$.   
Let  $Z$ be the finite-dimensional vector space, 
	$$Z=\R\oplus S(V^*) \oplus S^2(V^*) \oplus \cdots \oplus S^{d_0} (V^*).$$
We naturally identify $P_{\le - \lambda_1}(V)$ with a subspace $U\subset Z$, $U= \{ T\in Z: \varpi (T)\le -\lambda_1\}$,  
via polarization of homogeneous components.   
We equip each $S^k(V^*)$ and thus $Z$ and $U$ 
with the induced weight.  
If $V$ has an inner product, we equip each $S^k(V^*)$ and thus $Z$ and $U$ 
with the induced dual inner products.  

Let $\bfV$ denote the dual space,  $\bfV=U^*$.  
We equip $\bfV$ with the dual inner product induced from the inner product on $U$ and the dual weight.   
As the dual of a subspace $U\subset Z$,
  $\bfV$ is a quotient of $$Z^*=(\R^*\oplus S(V) \oplus S^2(V) \oplus \cdots \oplus S^{d_0} (V)).$$ 
  That is, $\bfV = Z^*/K$ where 
the kernel   $K$ consists of $v\in Z^*$ 
 such that $T(v)=0$ for all $T\in U.$ 
   Then $$K:=\{v\in \R^*\oplus S(V) \oplus S^2(V) \oplus \cdots \oplus S^{d_0} (V) :\varpi(v)<\lambda_1\}.$$
 
 Given any transversal $Y$ to $K$ in $Z^*$,
 the identification of $Y$ with $\bfV$  preserves weight by \cref{lem:compat}.  


We define  $\bfV'$ to be  the image of 
$$\{0\}\oplus S(V) \oplus S^2(V) \oplus \cdots \oplus S^{d_0} (V) $$
in $\bfV$ and hence is codimension-1.  Alternatively, in $Z$, consider the 1-dimensional subspace $C$ of constant functions in  $U=P_{\le - \lambda_1}(V)$.  Then $\bfV'$ coincides with the set of the elements of  $\bfV$  that vanish on $C$.


The subspace  $ P_{\le -\lambda_1}(V)$ of $Z$ contains the subspace  $V^*$ of linear functionals on $V$.  Then $V= S(V) =(V^*)^*$ is identified with a subspace of $\bfV$. Indeed, $\bfV$ is a quotient of $\R^*\oplus S(V) \oplus S^2(V) \oplus \cdots \oplus S^{d_0} (V)$ by $K$; we have that $K\cap S(V) = K\cap V = \{0\}$ and thus hence $V= S(V)$ injects $\bfV$.     
Let $\Pi_V^*\colon P_{\le - \lambda_1}(V)\to P_{\le - \lambda_1}(V)$ be the map $\Pi_V^*(f) = D_0f$.  We take $\Pi_V\colon \bfV\to \bfV$ to be the dual of $\Pi_V^*$.  Then the image of $\Pi_V$ is $S(V)=V$ in $\bfV$.

We take $\iota_V\colon V\to \bfV$ to be the evaluation map, $$\iota_V(v)(f ) = f (v).$$  We have that $\R^*$ is spanned by $\iota_V(0)$ (via the map $\phi \mapsto \phi(0)$) and \ref{LineNF1} follows.
Consider the constant function $1\colon V\to \R$.  Then the range of $\iota_V$ is contained in the codimension-1 affine subspace of $\bfV$ of vectors $\xi\in \bfV$ satisfying $\xi(1) = 1$.  This is precisely the space $\iota_V(0) + \bfV'$.  
We check that $\Pi_V\circ \iota_V$ is the identity map.  Indeed, if $f\in  P_{\le -\lambda_1}(V)$ and $v\in v$ then $$f(\Pi_V\circ \iota_V(v))=(\Pi_V^*f)(\iota_V(v)) = D_0f(\iota _V(v))= D_0f(v).$$
Since $ P_{\le -\lambda_1}(V)$ contains all linear functionals $V^*$, it follows that $\Pi_V\circ \iota_V(v)= v$.

Given $f \in \calP^{SR}(V,W)$, let $\bfL f \colon \bfV\to \bfW$ denote the adjoint of  $f ^*$.  

Properties \ref{LineNF2}, \ref{LineNF3}, \ref{LineNF5}, \ref{LineNF6}, \ref{LineNF7}, \ref{LineNF8}, and \ref{LineNF9} follow.

For  \ref{LineNF4}, it is clear that $f \mapsto \bfL f $  is continuous.  It follows from \ref{LineNF2} that the image of $\calG^{SR}(V)$ and $\calG^{SSR}(V)$ are closed subgroups of $\GL(\bfV)$; moreover  \ref{LineNF2} implies $\tau \colon f \mapsto \bfL f $ is an isomorphism.  
Let $D\tau\colon \Lie(\calG^{SSR}(V))\to \Lie ( \bfG^{SSR}(\bfV))$ and $D\tau\colon \Lie(\calG^{\ast}(V))\to \Lie ( \bfG^{\ast}(\bfV))$
 be the induced map of Lie algebras.  We have  $\Lie(\calG^{SSR}(V))$ and  $\Lie(\calG^{\ast}(V))$ are the vector spaces of strictly subresonant polynomial (resp. weight decreasing) maps $f \colon V\to V$; these are nilpotent lie algebras.  
 Consider $f\in \calG^{\ast}(V)$. 
 If $\phi\colon V\to \R$ has $\varpi(\phi)\le -\lambda_1$ satisfies $\phi\circ f = \lambda\phi$ we claim that $\lambda=0$. 
 Indeed, let $T\colon V^{\otimes k}$ be a homogeneous component of $\phi$.  
 Then
  $f^*$ has only zero eigenvalues. 
  It follows that $\Lie ( \bfG^{SSR}(\bfV))$ is nilpotent and relative to an appropriate basis are upper triangular and thus $ \bfG^{SSR}(\bfV)$ and $\calG^{\ast}(V) $ consist  of unipotent matrices.

\subsection{Block triangular form of the linearization}
A \emph{monomial function} on a vector space $W$ is a symmetric tensor $\phi \in S^k(W^*)$.   Given $f\in \calP^{SR}(W,V)$, $f^*$ and $\bfL f$ have a block structure using monomials as basis elements for $P_{\le -\lambda_1}(V)$ and $P_{\le -\lambda_1}(W)$.  
\begin{lemma}[Block triangular form of $f^*$]\label{lem:plplp|}
 Let $f\colon V\to W$ be a degree $d$ subresonant polynomial.  Let $\phi$ be a nonzero degree $k\ge 0$ monomial on $W$.  Then\begin{equation}\label{block}f^*\phi = \phi \circ (D_0 f \otimes \dots \otimes D_0 f) + \sum_{j=0}^{kd} \psi_j\end{equation}
where each $\psi_j$ is a degree $j$ homogenous polynomial such that 
\begin{enumerate}
\item $\varpi(\psi_j)\le \varpi(\phi)$ and 
\item if $\varpi(\psi_j)= \varpi(\phi)$ then $j>k$.
\end{enumerate}
\end{lemma}
\begin{proof}
When $k\ge 1$, the first conclusion follows from \cref{prop:subadd}.  If $\phi$ is degree $0$, we view $\phi\in \R$ and $f^*\phi =\phi.$

For the second conclusion, write $f$ as the sum of symmetric multilinear functions $f =F_0+  F_1 + \dots + F_d$ where $F_i = \wtd D^i f$.   
The polynomial function $f^*(\phi)\colon V\to \R$ is the sum    $$f^*(\phi)= \sum_{0\le i_1, \dots, i_k\le d}\psi_{i_1, \dots , i_k} $$ of monomials of the form
$$\psi_{i_1, \dots , i_k}= \phi\circ (F_{i_1}\otimes  \dots \otimes F_{i_k}).$$ 
Suppose $i_j= 0$ for some $1\le j\le k$.  By symmetry, we may assume $i_1 = 0$.   If $\mathbf{v}\in V^{\otimes (i_1 + \dots + i_k)}$ then 
\begin{align*}
\varpi (\psi_{i_1, \dots , i_k}(\mathbf{v}))&
\end{align*}
using that  each $F_{i_2}, \dots,  F_{i_j}$ is subresonant.  
It thus follows that $$\varpi(\psi_{i_1, \dots , i_k}) \le\varpi (\psi_{i_1, \dots , i_k}(\mathbf{v})) - \varpi( \mathbf{v})\le  \varpi(\phi) + \varpi( F_0)< \varpi(\phi)$$ since $\varpi( F_0)<0$ and $\varpi(\phi)>-\infty$.  

We conclude that  if $\varpi(\psi_{i_1, \dots , i_k}) = \varpi(\phi)$, then  $i_j\ge 1$ for every $1\le j\le k$.  
If $i_j = 1$ for every $j$ we obtain the term $\phi \circ (D_0 f \otimes \dots \otimes D_0 f)$ in \eqref{block}.
If at least one $i_j >1$ then $\psi_{i_1, \dots , i_k}= \phi\circ (F_{i_1}\otimes  \dots \otimes F_{i_k})$ has degree strictly larger than $k$ and the second conclusion follows.  
\end{proof}

The upper triangular form of $f^*$ and thus $\bfL f$ directly implies the unipotent elements of $\calG^{SR}(V) $ are of the form $\bfL(g)$ where the $\varpi(D_0g- \id)<0$.  In particular, we have the following.
\begin{corollary}\label{uniprad}
	$ \calG^{*}(V)$ is the unipotent radical of $ \calG^{SR}(V) $.
\end{corollary}

\subsection{Example}\label{sec:exam}
As an instructive example of the above construction, consider the {\SG} on $\R^3$ given by $$\varpi(e_1) = -3, \varpi(e_2) = -2,  \varpi(e_3) = -1$$ where $e_i$ are the standard basis vectors.  Every subresonant polynomial $f\colon \R^3\to \R^3$ is of the form
$$f(x,y,z) = (a_0+ a_1x + a_2 y + a_3 z+ a_4  yz + a_5 z^2 +a_6 z^3, b_0+  b_1 y+ b_2 z + b_3 z^2, c_0+ c_1z).$$

The subspace $P_{\le 3} (\R^3)$ has a basis the symmetrization of the following tensors
\begin{equation}\label{eq:basis1}\{1,e_1^*,e_2^*,e_3^*, e_2^*\otimes e_3^*, e_3^*\otimes e_3^*, e_3^*\otimes e_3^*\otimes e_3^*\}.\end{equation}
Relative to the ordered basis \eqref{eq:basis1}, $f^*$ has matrix $$
\left(
\begin{array}{ccccccc}
 1 &  {a_0} &  {b_0} &  {c_0} &  {b_0}  {c_0} &  {c_0}^2 &  {c_0}^3 \\
 0 &  {a_1} & 0 & 0 & 0 & 0 & 0 \\
 0 &  {a_2} &  {b_1} & 0 &  {b_1}  {c_0} & 0 & 0 \\
 0 &  {a_3} &  {b_2} &  {c_1} &  {b_0}  {c_1}+ {b_2}  {c_0} & 2  {c_0}  {c_1} & 3  {c_0}^2  {c_1} \\
 0 &  {a_4} & 0 & 0 &  {b_1}  {c_1} & 0 & 0 \\
 0 &  {a_5} &  {b_3} & 0 &  {b_2}  {c_1}+ {b_3}  {c_0} &  {c_1}^2 & 3  {c_0}  {c_1}^2 \\
 0 &  {a_6} & 0 & 0 &  {b_3}  {c_1} & 0 &  {c_1}^3 \\
\end{array}
\right)
$$
If we instead  order the basis  first by  weight and then by degree, we obtain the basis
\begin{equation}\label{eq:basis2}\{e_1^*,   e_2^*\otimes e_3^*,  e_3^*\otimes e_3^*\otimes e_3^*, e_2^*, e_3^*\otimes e_3^*, e_3^*, 1\}.\end{equation}
Relative to the basis \eqref{eq:basis2}, $f^*$ has the matrix 
\begin{equation}\label{eq:matrix}
\left(
\begin{array}{ccccccc}
  {a_1} & 0 & 0 & 0 & 0 & 0 & 0 \\
  {a_4} &  {b_1}  {c_1} & 0 & 0 & 0 & 0 & 0 \\
  {a_6} &  {b_3}  {c_1} &  {c_1}^3 & 0 & 0 & 0 & 0 \\
  {a_2} &  {b_1}  {c_0} & 0 &  {b_1} & 0 & 0 & 0 \\
  {a_5} &  {b_2}  {c_1}+ {b_3}  {c_0} & 3  {c_0}  {c_1}^2 &  {b_3} &  {c_1}^2 & 0 & 0 \\
  {a_3} &  {b_0}  {c_1}+ {b_2}  {c_0} & 3  {c_0}^2  {c_1} &  {b_2} & 2  {c_0}  {c_1} &  {c_1} & 0 \\
  {a_0} &  {b_0}  {c_0} &  {c_0}^3 &  {b_0} &  {c_0}^2 &  {c_0} & 1 \\
\end{array}
\right).\end{equation}
Observe that diagonal entries of are all induced by the linear part of $f$ and the off-diagonal terms either decrease weight or increase degree if they preserve weight.  

Formally, $\bigl(P_{\le 3} (\R^3)\bigr)^*$ is quotient of $$\R^*\oplus (\R^3)^*\oplus S^2((\R^3)^*)\oplus S^3((\R^3)^*)$$
which we can identify with the subspace spanned by the symmetrization of the following tensors
\begin{equation}\label{eq:basis3}\{e_1, e_2\otimes e_3,e_3\otimes e_3\otimes e_3, e_2, e_3\otimes e_3, e_3, 1^*\}.\end{equation}
Relative to the basis \eqref{eq:basis3}, $\bfL f$  is the transpose of the \eqref{eq:matrix}.  In particular, in blocks of the same weight, off diagonal terms of $\bfL f$ decrease degree.

\section{Subresonant structures}\label{sec4}
Fix $\lambda_1<\dots<\lambda_\ell<0$ 
and  $r> \lambda_1/\lambda_\ell\ge 1$.

\subsection{Subresonant structures on manifolds}
Let $N$ be a $C^r$ manifold; as it suffices in all applications we are concerned with, we will further assume that $N$ is diffeomorphic to $\R^n$ for some $n$.
\begin{definition}
	A \emph{complete $C^r$ subresonant structure} on $N$ with weights   $\lambda_1<\dots<\lambda_\ell<0$ is the following data: 
	\begin{enumerate}
		\item a family of compatible {\SGs} $\varpi_x$ taking values $\lambda_1<\dots<\lambda_\ell<0$ defined on  $T_xN$ for every $x\in N$, and 
		\item an atlas $\calA_x= \{h\colon T_xN\to N\}$ of $C^r$ diffeomorphisms defined for every $x\in N$
	\end{enumerate}
	such that 
	\begin{enumerate}[label=(\alph*),ref=(\alph*)]
		\item the group $\calG^{SR}(T_xN)$ acts transitively on each $\calA_x$ by precomposition, and 
		\item for all $x,y\in N$, every $h_x\in \calA_x$, and every $h_y\in \calA_y$, the transition map $h_x\inv\circ h_y$ is an element of $\calP^{SR}(T_yN, T_xN)$. 
	\end{enumerate}
\end{definition}
Observe that if $N$ has a complete $C^r$ subresonant structure, then the parameterized family of flags $x\mapsto \calV_{\omega_x}$  is $C^{r-1}$.

\subsection{Subresonant diffeomorphisms}
Let $N_1$ and $N_2$ be $C^r$ manifolds equipped with complete subresonant structures with weights in $\lambda_1<\dots<\lambda_\ell<0$.

We say the subresonant structures are \emph{compatible} if the multiplicities of $\varpi_x$ and $\varpi_y$ coincided for some (and hence all) $x\in N_1$ and $y\in N_2$.  
\begin{definition}
	Suppose $N_1$ and $N_2$ have compatible subresonant structures.  
	A $C^r$ diffeomorphism  $f\colon N_1\to N_2$ is \emph{subresonant} if for any (and hence all) $x\in N_1$, $y\in N_2$, $h_x\in \calA_x$, and $h_y\in \calA_y$,
	$h_y\inv \circ f \circ h_x$
	is an element of $\calP^{SR}(T_xN_1, T_yN_2)$.
\end{definition}

\subsection{Linearization of subresonant diffeomorphisms}
Let $N_1$ and $N_2$ be $C^r$ manifolds equipped with complete subresonant structures taking weights  $\lambda_1<\dots<\lambda_\ell<0$.  
Fix $x\in N_1$ and $h_x\in \calA_x$.
Let $$V_1 := \{\phi\circ h_x\inv: \phi \in P_{\le -\lambda_1} (T_xN_1)\}.$$
As a vector space of real-valued functions on $N_1$, note that  $V_1$ is independent of the choice of $x\in N_1$ and $h_x\in \calA_x$.  
Similarly define $V_2:= \{\phi\circ h_y\inv: \phi \in P_{\le -\lambda_1}( T_yN_2)\}$ for any choice of $y\in N_2$ and $h_y\in \calA_y$.

Let $f\colon N_1\to N_2$ be a subresonant $C^r$ diffeomorphism.  
We have $f^*\colon \phi\mapsto \phi\circ f$ is a linear map from $f^*\colon V_2\to V_1$

Let $\bfL N_1 = V_1^*$ and $\bfL N_2 = V_2^*$.  
We embed $N_i$ into $\bfL N_i$,  
$\iota_i\colon N_i \to \bfL N_i$, by
$$\iota_i(x) (\phi) = \phi(x).$$
The adjoint of $f^*$ is then a linear map $\bfL f\colon \bfL N_i\to \bfL N_2$.
Moreover, we have
$$\bfL f \circ \iota_1= \iota_2\circ f.$$

Note that $T_xN_1$ has a subresonant structure relative to which each $h_x\in \calA_x$ is a subresonant diffeomorphism.  
Fix $x\in N_1$, $y\in N_2$, $h_x\in \calA_x$, and $h_y\in \calA_y$.  
Then $$\bfL(h_y\inv \circ f\circ h_x),$$ the linearization of the subresonant polynomial $h_y\inv \circ f\circ h_x$ as defined in \cref{thm:LinSRPoly}, 
coincides with $$\bfL(h_y\inv) \circ \bfL f\circ \bfL h_x.$$

In particular, while the vector spaces $\bfL N_1$ and $\bfL N_2$ and the linear map $\bfL f$ are intrinsically defined independent of the choice of $x\in N_1$ or $y\in N_2$, we will often instead study the map induced by a choice of $x\in N_1$, $y\in N_2$, $h_x\in \calA_x$, and $h_y\in \calA_y$.  

\subsection{Induced inner products}
For each $x\in N_1$, we have a natural identification $$P_{\le -\lambda_1} (T_xN_1) = \R^*\oplus T^*_xN_1 \oplus 
S^2(T^*_xN_1) \oplus \dots \oplus S^k(T^*_xN_1 ).$$
Given $h_x\in \calA_x$, the map $\phi\mapsto \phi\circ h_x\inv$ induces an isomorphism of vector spaces $P_{\le -\lambda_1} (T_xN_1) \to V_1.$
An inner product on $T_xN_1$ induces an inner product on each $S^j(T^*_xN_1 )$ which thus induces an inner product on $V_1$ and on $\bfL N_1 = V_1^*$; however this inner product depends both on the choice of $x\in N_1$ and on the choice of $h_x\in \calA_x$.

\section {Lyapunov exponents and forwards regularity} \label{sec5}
Let $\{W^i\}_{i\in \Z}$ be a sequence of  finite-dimensional inner product spaces of constant dimension $d= \dim W^i$.  Consider a sequence of invertible linear maps $A_i\colon W^i\to W^{i+1}$ and for $n\ge 1$, write $$A^{(n)}_i:= A_{n+i-1} \circ \dots \circ A_i.$$
We will always assume the sequence satisfies \begin{equation} \label{eq:controlednorm} \limsup_{n\to \infty} \frac 1 n \log \|A_n \| = \limsup_{n\to \infty} \frac 1 n \log \|A_n \inv \| = 0.\end{equation}
Given $v\in W^0$, set \begin{equation}\label{eq:lyapSG}\varpi(v) = \limsup_{n\to \infty}\frac 1 n \log \| A^{(n)}_0 v\|.\end{equation}
By \eqref{eq:controlednorm}, we have have $-\infty< \varpi(v) <\infty$ for all $v\neq 0$.  
Moreover, it is clear that $\varpi$ defines a {\SG} on $W^0$ (c.f. \cite[Proposition 1.3.1]{MR2348606}), which we refer to as the \emph{Lyapunov} {\SG}.  Let $\lambda_1< \dots< \lambda_\ell$ denote the image of $W^0\sm \{0\}$ under  $\varpi$ and let $m_i$ be the associated multiplicities.  The numbers $\{\lambda_i\}$ are called the \emph{Lyapunov exponents} of the sequence $\{A_i\}$.  
\begin{definition}
The sequence of linear maps $\{A^{(n)}_0\}_{n\in \N}$  is \emph{forwards regular} if the limit $$\lim_{n\to \infty } \frac 1 n \log |\det A^{(n)}_0|$$ exists (and is finite) and is equal to $\sum_{i=1}^\ell m_i \lambda_i.$
\end{definition}
We remark if $\{A^{(n)}_0\}_{n\in \N}$ is forwards regular then $\{A^{(n)}_i\}_{n\in \N}$ is forwards regular for any $i\in \Z$ as are all exterior powers of $\{A^{(n)}_0\}_{n\in \N}$.   

\subsection{Properties of forwards regular sequences}
Let $\{A_i\}_{n\in \N}$ be a sequence of invertible linear transformations $A_i\colon W^i\to W^{i+1}$ as above.  Let $\lambda_1<\dots<\lambda_\ell$ be the Lyapunov exponents with associated multiplicities $m_1, \dots, m_\ell$.  
Let $$\calV^0 = \{V_1\subset V_2\subset \dots \subset V_\ell\}$$ be the associated filtration of $W^0$. 
As in \cref{def:adapbasis}, an ordered basis  $\{e_i\}$ of $W^0$ is  \emph{adapted} to the flag $\calV^0$ if for each $1\le j\le \ell$,
$$V_j = \Span\{e_i: 1\le i \le {m_1 + \dots + m_j} \}.$$
A direct sum decomposition  $W^0= \bigoplus E_j$ of $W^0$ is adapted to the flag $\calV_0$ if for each $1\le j \le \ell$, 
$ V_{j-1} \cap E_j=\{0\}$ and 
$$V_j = V_{j-1} \oplus E_j.$$
For instance, given an inner product on $W^0$ we may take $E_j= V_{j-1}^\perp \cap V_j$ or $$E_j = \Span\{e_i: {m_1 + \dots + m_{j-1}}+1 \le i \le {m_1 + \dots + m_j} \}$$ for any basis $\{e_i\}$ adapted to the flag $\calV^0$.  

We have the following well-known facts about forward regular sequences.
\begin{proposition}
Let $\{A^{(n)}_0\}_{n\in \N}$ be forwards regular. 
\begin{enumerate}[label=(\arabic*),ref=(\arabic*)]
\item For every splitting $W^0= \bigoplus E_j$  adapted to the flag $\calV^0$ and every $1\le j\le \ell$, 
$$\lim_{n\to \infty} \frac 1 n \log \|{A^{(n)}_0} v \| = \lambda_j$$
uniformly over $v$ in the unit sphere in $E_j$.
\item For every basis $\{e_i: 1\le i \le d\}$  of $W^0$ adapted to the flag $\calV^0$ and any disjoint subsets $I,J \subset \{1, 2, \dots ,d \}$,
write $E_I =  \Span \{e_i: i\in I\}$ and $E_J =  \Span \{e_i: i\in J\}$.  Then 
$$\lim_{n\to \infty}\frac 1 n \log \sin \angle\bigl({A^{(n)}_0} E_I, {A^{(n)}_0} E_J\bigr) =0.$$

\end{enumerate}
\end{proposition}

\begin{corollary}\label{cor:LyapControl}
Suppose the sequence $\{A^{(n)}_0\}_{n\in \N}$ is forwards regular.
Fix a direct sum decomposition  $W^0=\bigoplus E_j$ and a basis $\{e_i\}$  adapted to the flag $\calV^0$.

For any $\epsilon>0$, there is $C_\epsilon>1$ (depending on the choice of $\bigoplus E_j$ and $\{e_i\}$) so that the following hold for all $n\ge 0$ and $k\ge 0$:
\begin{enumerate}[label=(\arabic*),ref=(\arabic*)]
\item  For every $v\in E_j$,
\begin{equation}\label{growth controls}
\frac 1 {C_\epsilon} e^{n(\lambda_j -\epsilon)} \|v\| \le \|{A^{(n)}_0} v \| \le  {C_\epsilon} e^{n(\lambda_j +\epsilon)} \|v\|. 
\end{equation}
\item  If $\wtd v_k = {A^{(k)}_0} v$ then 
$$\frac 1 {{C_\epsilon}^2} e^{-2\epsilon k} e^{n(\lambda_j -\epsilon)} \|\wtd v_k \| \le \|{A^{(n)}_k} \wtd v_k \| \le  {C_\epsilon} ^2 e^{2\epsilon k} e^{n(\lambda_j +\epsilon)}\|\wtd v_k \|. 
$$
\item\label{LyapCont3} Let $d= \dim W^0$.  If $v= \sum_{i=1}^d v_i e_i$, then  
 $$ \frac 1 {C_\epsilon} e^{-\epsilon n} \max \bigl \{ v_i  \|{A^{(n)}_0} e_i \|\bigr\} \le \|{A^{(n)}_0} v \|\le d  \max\bigl \{ v_i  \|{A^{(n)}_0} e_i \|\bigr\}.$$
\end{enumerate}
\end{corollary}
 
\subsection{One-sided Lyapunov metric}
Let $\{A_i\}_{n\in \N}$ be a sequence of invertible linear transformations $A_i\colon W^i\to W^{i+1}$ as above.  Let $\lambda_1<\dots<\lambda_\ell$ be the Lyapunov exponents with associated multiplicities $m_1, \dots, m_\ell$.  
Let $$\calV^0 = \{V_1\subset V_2\subset \dots \subset V_\ell\}$$ be the associated filtration of $W^0$. 

Let $d= \dim W^i$.  
\begin{proposition}\label{absLyapu}
Fix $\epsilon>0$.   Suppose $\{A_i^{(n)}\}_{n\in \N}$ is forwards regular.   There  exists a constant $L_\epsilon$ and, for every $k\in \{0,1,2, \dots\}$, an inner product $\iprod{\cdot,\cdot}'_k$ with associated norm $\|\cdot\|'_k$ on $W^k$ with the following properties: for every $v\in W^0$,
\begin{enumerate}
\item \label{cock1}$\|A_0^{(n)}v\|_n'\le e^{n (\varpi(v)+\epsilon)}\|v\|_0'$  for every $n$;
\item\label{cock2} $  \|v\| \le  \|v\|'_n \le L_\epsilon e^{\frac {\epsilon }{2} n} \|v\|$ for every $n$;
\item  \label{cock3}$\|A_0^{(n)}v\|_n'\ge e^{n (\varpi(v)-\epsilon)}\|v\|_0'$  for all sufficiently large (depending on $v$) $n$.  
\end{enumerate}
\end{proposition}
\begin{proof}
Given $v,w\in W^k$, set 
	$$\iprod{v,w}'_k:=\sum_{n=0}^\infty e^{-n(2\epsilon +\varpi(v) + \varpi(w))} \iprod {A_k^{(n)}v,A_k^{(n)}w}.$$
From \cref{cor:LyapControl}, this expression converges and the lower bound in \eqref{cock2} holds.
For the upper bound in  \eqref{cock2}, we apply  use  $\epsilon/2$ in the place of $\epsilon$ in \cref{cor:LyapControl} and obtain
\begin{align*}
\left(\| A_0^{(k)} v\|'_k\right)^2 &=
\sum_{n=0}^\infty e^{-n(2\epsilon +2\varpi(v) )} \| A_k^{(n)}v\|^2\\
&\le \sum_{n=0}^\infty e^{-n(2\epsilon +2\varpi(v) )} C_{\frac \epsilon 2} ^2 e^{\epsilon k}e^{2n(\varpi(v) +\frac \epsilon 2)}\| v \|^2 \\
&\le C_{\frac \epsilon 2} ^2 e^{\epsilon k} \sum_{n=0}^\infty e^{-n\epsilon }\| v \|^2.
\end{align*}
Set $L_\epsilon :=  \left(C_{\frac \epsilon 2} \sum_{n=0}^\infty e^{-n\epsilon }\right)^{\frac 1 2}.$
\eqref{cock3} follows from  \cref{cor:LyapControl} 
and the upper bound in \eqref{cock2}.  

For  \eqref{cock1}, given $v\in W^0$ we have 
\begin{align*}
\left(\| A_0^{(k)} v\|'_k\right)^2 &=
\sum_{n=0}^\infty e^{-n(2\epsilon +2\varpi(v) )} \| A_k^{(n)}v\|^2\\
&\le  \sum_{n=0}^\infty e^{-(n-k)(2\epsilon +2\varpi(v) } \| A_0^{(n)}v\|^2\\
&=  e^{k(2\epsilon +2\varpi(v) )}\sum_{n=0}^\infty e^{-n(2\epsilon +2\varpi(v) } \| A_0^{(n)}v\|^2\\
&=  e^{k(2\epsilon +2\varpi(v) )} \left(\|   v\|'_0\right)^2
\end{align*}

\end{proof}

\subsection{Regularity and Lyapunov exponents for block triangular matrices}
We now assume each $W^i$ is equipped with an ordered orthonormal basis adapted to a Lyapunov filtration.  In this way, we naturally identify each $W^i$ with $\R^d$ equipped with the standard inner product.  
This simplifies the statement and proof of the following  estimate.
\begin{lemma}\label{unipotentisok}
Let $V$ be a subspace of $\R^d$ and let $L_i\colon \R^d\to \R^d$ be a sequence of invertible matrices with $L_i(V)= V$ for every $i$.  
Let $W= V^\perp$ and write each $L_i$ in block form 
$$L_i = \left(\begin{array}{cc}A_i & U_i \\0 & B_i\end{array}\right)$$
relative to $V\oplus W$ so that $A_i\colon V\to V$, $B_i\colon W\to W$ and $U_i\colon W\to V$ for every $i$.  

Suppose the sequences $\{A_{0}^{(n)}\}$ and  $\{B_{0}^{(n)}\}$ are forwards regular and that 
 $$\lim_{n\to \infty } \frac 1 n \log ^+\|U_{n}\| =0.$$
 Then the sequence $\{L_{0}^{(n)}\}$ is forwards regular and the Lyapunov exponents of $\{L_{0}^{(n)}\}$ counted with multiplicity coincide with the union of the Lyapunov exponents of the sequences $\{A_{0}^{(n)}\}$ and $\{B_{0}^{(n)}\}$ counted with multiplicity.  
\end{lemma}

\begin{proof}
Let $\{\wtd L_i\}$ be the sequence of matrices with block form $$L_i = \left(\begin{array}{cc}A_i & 0 \\0 & B_i\end{array}\right).$$
Then $\{\wtd L_{0}^{(n)}\}$ is forwards regular and $$\lim_{n\to \infty } \frac 1 n \log |\det L^{(n)}_0|= 
\lim_{n\to \infty } \frac 1 n \log |\det\wtd  L^{(n)}_0|.$$
It thus suffices to show the sequences $\{\wtd L_{0}^{(n)}\}$ and $\{L_{0}^{(n)}\}$ have the same Lyapunov exponents.  

Let $\lambda_1\ge \lambda_2\ge \dots \ge \lambda_p$ and $\eta_1\ge \eta_2\ge \dots \ge \eta_q$, $p+q= d$,  denote the Lyapunov exponents of the sequences $\{A_{0}^{(n)}\}$ and $\{B_{0}^{(n)}\}$, respectively, listed with multiplicities.  Let $\{e_i\}$ and $\{f_i\}$ be bases of $V$ and $W$, respectively, adapted to the Lyapunov {\SGs}; in particular, we have 
\begin{equation}\label{growthindir}\lim_{n\to \infty} \frac 1 n \log\|A_0^{(n)} e_i \| = \lambda_i\quad \quad 
\lim_{n\to \infty} \frac 1 n \log\|B_0^{(n)} f_j \| = \eta_j.\end{equation}
Given $n\ge 0$, write  
$$e_i^n = \frac {A_0^{(n-1)}e_i}{\|A_0^{(n-1)}e_i\|}, \quad \quad 
f_{i}^n= \frac {B_0^{(n-1)}f_j}{\|B_0^{(n-1)}f_j\|}$$
for the renormalized image bases.
Given $1\le j \le q$ and $1\le i\le p$, let $u_{n,i,j}\in \R$ be such that $$U_n(f_j^n) = \sum u_{i,j,n}e_i^n.$$

Combining \cref{cor:LyapControl}\ref{LyapCont3} with the fact that  $\lim_{n\to \infty }\frac 1 n \log ^+\|U_{n}\| =0$, we have 
\begin{equation}\label{eq:slowunip}\lim_{n\to \infty } \frac 1 n \log  \max \{ |u_{i,j,n}|: 1\le i \le p, 1\le j\le q\} = 0.\end{equation}

Fix $1\le j\le q$.  Consider  any $1\le i \le p $ for which $\eta_j< \lambda_i$.
 For such  $j$ and $i$, write 
\begin{equation}
c_{i,j} = \sum _{k= 1}^\infty  \| A_0^{k} e_{i}\|\inv  u_{i,j,k-1} \|B_0^{(k-1)} f_j\|.
\end{equation}
Fix $\epsilon>0$.  By \cref{cor:LyapControl} and   \eqref{eq:slowunip}, there is $C_\epsilon>1$ such that 
\begin{align*}
|c_{i,j}|&\le 
\sum _{k= 1}^\infty  \| A_0^{k} e_{i}\|\inv  |u_{i,j,k-1}| \|B_0^{k-1} f_j\|\\
&\le 
\sum _{k= 1}^\infty  \big( C_\epsilon e^{-k(\lambda_i -\epsilon)}\bigr) \bigl(C_\epsilon  e^{k \epsilon }\bigr) \bigl(C_\epsilon e^{(k-1)(\eta_j +\epsilon)}\bigr) 
\end{align*}
which converges assuming $\lambda_i>\eta_j$ and $\epsilon$ is taken  sufficiently small so that $\epsilon <(\lambda_i-\eta_j)/3$.

Let 
 $$\wtd f_j = f_j - \sum_{\lambda_i >\eta_j} c_{i,j} e_i.$$
We claim \begin{equation}\label{eq:correctgrowth}\lim_{n\to \infty} \frac 1 n \log\|L_0^{(n)}  \wtd f_j \| = \eta_j.\end{equation}
Indeed, for every $1\le i\le p$ such that  $\lambda_i\le \eta_j$, the $e_i^n$ component of 
$L_0^{(n)}  \wtd f_j$
is 
\begin{align*}
\sum_{k=1}^{n} u_{i,j,k-1} \|B^{(k-1)}_0 f_j\| A_k^{(n-k)} (e_{i}^k)
\end{align*}
which, for every $n\ge 1$, has norm bounded above by 
\begin{align}
\sum_{k=1}^{n} C_\epsilon e^{\epsilon k} & C_\epsilon e^{(k-1)(\eta _j+\epsilon)} e^{2k\epsilon } {{C_\epsilon} ^2}
e^{(n-k)(\lambda_i+\epsilon)}
\notag \\
 &\le 
\sum_{k=1}^{n} {C_\epsilon}^4  e^{3 \epsilon k} e^{(n-1)(\eta _j+\epsilon)} \notag \\
 &\le 
\frac {e^{3\epsilon n}- 1}{ 1 -e^{-3\epsilon } } {C_\epsilon}^4  e^{(n-1)(\eta _j+\epsilon)}.  \label{growthslowdir}
\end{align}
On the other hand, for $1\le i\le p$ for which $\lambda_i>  \eta_j$, the $e_i^n$ component of  $L_0^{(n)}  \wtd f_j$ is 
\begin{align*}
&\sum_{k=1}^{n} u_{i,j,k-1}  \|B^{(k-1)}_0 f_j\| A_k^{n-k} (e_{i}^k) 
-A_0^{(n)}c_{i,j} e_i\\
&=\sum_{k=1}^{n} u_{i,j,k-1}  \|B^{(k-1)}_0 f_j\| A_k^{n-k} (e_{i}^k) \\
& \quad \quad \quad  -
 \sum _{k= 1}^n \| A_0^{k} e_{i}\|\inv  u_{i,j,k-1} \|B_0^{(k-1)} f_j\| A_0^{n}   (e_i)\\
& \quad \quad \quad
-\sum _{k= n+1}^{\infty}  \| A_0^{k} e_{i}\|\inv  u_{i,j,k-1} \|B_0^{(k-1)} f_j\|  A_0^{n} (e_i)\\
=& \sum_{k=1}^{n} u_{i,j,k-1}  \|B^{(k-1)}_0 f_j\| A_k^{n-k} (e_{i}^k)\\
& \quad \quad \quad
-\sum _{k= 1}^n   u_{i,j,k-1} \|B_0^{(k-1)} f_j\|  A_k^{n-k}   (e_i^k)\\
& \quad \quad \quad
-\sum _{k= n+1}^{\infty} \| A_0^{k} e_{i}\|\inv  u_{i,j,k-1} \|B_0^{(k-1)} f_j\|  A_0^{n} (e_i)\\
=& 
-\sum _{k= n+1}^{\infty} u_{i,j,k-1} \|B_0^{(k-1)} f_j\|  \| B_n^{k} (e_i^n)\| \inv e_i^n
\end{align*}
which has norm bounded above by 
\begin{align}\sum_{k=n+1}^{\infty } &C_\epsilon e^{\epsilon k} C_\epsilon e^{(k-1)(\eta _j+\epsilon)} {C_\epsilon ^2} e^{2\epsilon n}
e^{(k-n)(-\lambda_i+\epsilon)} \notag \\
&\le \sum_{k'=0}^{\infty } C_\epsilon e^{\epsilon (k' +n+1)}  C_\epsilon e^{(k'+n)(\eta _j+\epsilon)} {C_\epsilon ^2}e^{2\epsilon n}
e^{(k'+1)(-\lambda_i+\epsilon)}\notag \\
&\le \sum_{k'=0}^{\infty } {C_\epsilon}^4 e^{\epsilon (3n+1)} e^{n(\eta _j+\epsilon)}      e^{-\lambda_i+\epsilon} e^{k'(\eta _j-\lambda_i +3\epsilon)} \notag  \\
& \le  e^{n\eta _j+(4n+1)\epsilon} C \label{growthfastdir}
 \end{align}
 where, having taken $\epsilon>0$ sufficiently small,   $C<\infty$ is independent of $n$.  
 Since the $f_i^n$ component of  $L_0^{(n)}  \wtd f_j$ is $ B_0^{(n)} f_j$, the claim in \eqref{eq:correctgrowth} then follows from \eqref{growthindir},  estimates \eqref{growthslowdir} and \eqref{growthfastdir}, and   the arbitrariness of $\epsilon>0$  combined with \cref{cor:LyapControl}\ref{LyapCont3}.
\end{proof}

\subsection{Regularity and Lyapunov exponents for linearization of subresonant polynomial maps.}
Let $\{V^i\}_{i\in \Z}$ be a sequence of isomorphic, finite-dimensional  inner product spaces equipped with compatible {\SGs} $\varpi_i$.
We suppose each {\SG} takes negative values $\lambda_1<\dots <\lambda_\ell<0.$
Recall we write $\bfV^i= \bigl(P_{\le -\lambda_1}(V^i)\bigr)^*$

We equip each $P_{\le -\lambda_1}(V^i)$ with a basis of monomials.  We order this basis first relative to weight $\varpi$, and then relative to degree in each space of the same weight.  
Given $f\colon \calP^{SR}(V^i, V^{i+1})$, \cref{lem:plplp|} implies $f^* \colon P_{\le -\lambda_1}(V^{i+1}) \to P_{\le -\lambda_1}(V^i) $ has a block triangular structure.  It follows from \cref{lem:plplp|}  that the adjoint $\bfL f\colon \bfV^i\to\bfV^{i+1}$ also has a block triangular structure relative to which the diagonal blocks are the maps induced by tensor powers of $D_0f$; moreover the off-diagonal blocks have norm dominated by $\|f\|_P$ where $\|f\|_P$ is as in \eqref{eq:normofpoly}.

Applying \cref{unipotentisok} recursively, we immediately obtain the following. 

\begin{corollary}
For each $i\in \Z$, let $\{f_i\colon V^i\to V^{i+1}\}$ be an invertible subresonant polynomial map.  
Let $A_i= D_0f_i$ and let $\bfA_i = \bfL f_i$.
Suppose that 
\begin{enumerate}
\item $\limsup _{n\to \infty}\frac 1 n \log  \|f_i\|_P= 0$;
\item the sequence $\{A_0^{(n)} \}$ is forwards regular;
\item $\lim \frac 1 n \log \| A_0^{(n)} v\| = \varpi_0(v)$ for every $v\in V^0$.  
\end{enumerate}

Then the sequence of linear maps $\{\bfA_0^{(n)}\}$ is forwards regular and the Lyapunov exponents of the sequence 
$\{\bfA_0^{(n)}\}$ are all expressions of the form
$$ \lambda_1 \le  \sum_{i = 1} ^\ell n_i \lambda_i \le 0 $$
where $ n_i$ are non-negative integers.

Moreover, the exponent $0$ has multiplicity 1 and the Lyapunov exponents for the restriction of $\{\bfA_0^{(n)}\}$ to the codimension-1 subspace $(\bfV^0)'$
are all expressions of the form
$$ \lambda_1 \le  \sum_{i = 1} ^\ell n_i \lambda_i < 0 $$
where $ n_i$ are non-negative integers that are not all identically zero.  
\end{corollary}

\section{Temperedness}\label{sec6}
One minor nuisance in the formulation of results in  \cite{MR3642250} involves various notions of slow-exponential growth of certain estimates along orbits of the dynamics.  Specifically, in  \cite{MR3642250}, the $C^{r}$-norm of the localized dynamics is assumed $\epsilon$-tempered (see \cref{def:e-temp} below); the normal form change of coordinates and the size of the polynomial dynamics in these normal form coordinates are then shown to be $\kappa$-tempered where $\kappa$ is  depends on $r$ and $\epsilon$.  

However, an elementary argument shows that given any $M>0$, any $M$-tempered function is automatically $\epsilon$-tempered for all $\epsilon>0$; see \cref{lem:alltempsame}.   This allows us to avoid quoting precise temperedness estimates from \cite{MR3642250}; see \cref{cor:temp}. 
\subsection{Temperedness of functions}
Let $(X,\mu)$ be a probability space and let $f\colon X\to X$ be a measurable, $\mu$-preserving transformation.  
\begin{definition}\label{def:e-temp}
Let $\phi\colon X\to [1,\infty)$ be a measurable function.  Given $\epsilon>0$, we say $\phi$ is \emph{$\epsilon$-tempered} or 
\emph{$\epsilon$-slowly growing} if for almost every $x\in X$,
$$\sup \{ e^{-\epsilon n} \phi (f^n(x)) :n \ge0  \}<\infty.$$
We say that $\phi$ is \emph{tempered} 
if it is $\epsilon$-tempered for all $\epsilon>0$.  
\end{definition}
We observe that if $\phi$ is $\epsilon$-tempered then 
$\limsup_{n\to\infty} \frac 1 n \log \phi(f^n(x)) \le  \epsilon.$
Also, $\phi$ is tempered if and only if 
$\limsup_{n\to\infty} \frac 1 n \log \phi(f^n(x)) =0.$

We have the following standard construction which often is useful to provide dynamical bounds of tempered objects: Suppose $\phi\colon X\to [1,\infty)$ is an $\epsilon$-tempered, measurable function.  
Take \begin{equation}\label{eq:tempC} C_\epsilon(x):= \sup \{ e^{-\epsilon n} \phi (f^n(x)) :n  \ge 0 \}.\end{equation}
 Then $C _\epsilon(x)<\infty$ a.s., $C _\epsilon\colon X\to \R$ is measurable, and for almost every $x\in X$, we have $\phi(x)\le C _\epsilon(x)$ and 
\begin{equation}\label{eq:tempgrowth}
	C _\epsilon(f(x))\le e^\epsilon C_\epsilon(x).
\end{equation}
In particular, 
$$\phi(f^n(x))\le C_\epsilon(f^n(x))\le e^{n\epsilon} C_\epsilon (x).$$

\subsection{At-most-exponential growth implies temperedness}
The following lemma implies that temperedness of a function follows once $\phi$ grows at most exponentially along orbits.  
\begin{lemma}\label{lem:lllkop}\label{lem:alltempsame}
Let $f\colon X\to X$ be $\mu$-preserving.  
Let $\phi\colon X\to [1,\infty)$ be a function such that for some $M\ge 1$, $$\phi(f^n(x)) \le M^n \phi(x)$$
for almost every $x$.  
Then $\phi$ is  tempered.  
\end{lemma}
\begin{proof}
Fix $\ell_0>0$ such that, writing $$Y= Y_{\ell_0}:= \{x: \phi(x) \le \ell_0\}>0,$$ we have $\mu(Y)>0$.  Given $y\in Y$, let $n(y) = \min\{j\ge 1: f^j(y) \in Y\}$ be the first return time to $Y$.  
Then \begin{equation}\label{eq:ppp}\int_Y n(y) \ d \mu(y) = 1<\infty.\end{equation}
Let $g\colon Y\to Y$ denote the first return map to $Y$: $g(y) = f^{n(y)}(y)$.
From \eqref{eq:ppp} and the pointwise ergodic theorem (for the dynamics of $g\colon Y\to Y$),  we have \begin{align*}
 	\frac{1}{j}\log \left(M^{n(g^j(y))} \ell_0 \right) = 
 	\frac{1}{j} \bigl(\log(\ell_0)+ (\log M)  n(g^j(y))\bigr)
 	 \to 0
\end{align*}
 as $j\to \infty$ for almost every  $y\in Y$.

Fix $x\in X$.  By Poincar\'e recurrence, for almost every $x\in X$, there exists $\ell_0>0$ such  that $x\in Y_{\ell_0}$ and $f^{n_j}(x)\in Y_{\ell_0}$ for infinitely many $n_j\in \N$.
Fix such $x\in X$ and $\ell_0$ and let $$0= n_0< n_1< n_2< \dots$$denote the subsequent times for which $f^{n_i}(x) \in Y_{\ell_0}$.  
Then for $n_j\le n < n_{j+1}$ we have $j\le n_j\le n$ and 
\begin{align*}
	\frac 1 n \log \phi (f^n(x)) 
&\le \frac {1}{n} \log \left(M^{n-n_j}\ell_0 \right )
\le \frac {1}{n} \log \left(M^{n_{j+1}-n_j}\ell_0 \right )
\\
	&= \frac {1}{n} \log  \left(M^{n(f^{n_j}(x))}\ell_0 \right )
	\le \frac {1}{j} \log  \left(M^{n(f^{n_j}(x))}\ell_0\right)
		\\
		&= \frac {1}{j} \log \left(M^{n(g^j(x))}\ell_0\right)
	\end{align*}
whence \[\frac 1 n \log  \phi (f^n(x))\to 0.\qedhere\]
\end{proof}

From \eqref{eq:tempgrowth} and \cref{lem:alltempsame} thus obtain the following
\begin{corollary}\label{cor:temp}
	Let $f\colon X\to X$ be $\mu$-preserving and let $\phi\colon X\to [1,\infty)$ be a function.
	Then $\phi$ is tempered if and only if $\phi$ is $\epsilon$-tempered for some $\epsilon>0$.  
\end{corollary}

Note that while $\epsilon$-temperedness implies temperedness, the function  $C_\epsilon(x)$ in \eqref{eq:tempC} will in general depend on the choice of $\epsilon>0.$

\section{Contracting foliations and subresonant structures}\label{sec7}
Let $M$ be a (possibly non-comapact) manifold equipped with continuous Riemannian metric.    Fix $r> 1$, let $f\colon M\to M$ be a $C^r$ diffeomorphism.  Later, we will impose additional criteria on $r$ depending on the Lyapunov exponents of a certain cocycle.  

  Let $\mu$ be an ergodic, $f$-invariant Borel probability measure on $M$.  
We assume that relative to the metric and measure  that \begin{equation}\label{eq:integrability} \int \log \|D_x f\| \ d \mu(x) <\infty, \quad \quad \int \log \|D_x f\inv\| \ d \mu(x) <\infty.\end{equation}

 \subsection{Contracting foliations}
\begin{definition}\label{def:contfol}
Let $\calF$ be a (possibly non-measurable) partition of $M$.  We say $\calF$ is a \emph{contracting, $C^r$-tame, $f$-invariant, measurable foliation} if the following hold for almost every $x\in M$.
\begin{enumerate}
\item The atom $\calF(x)$ of $\calF$ is an injectively immersed, $C^r$ submanifold of $M$ which is diffeomorphic to $\R^d$.
 \item  $f(\calF(x)) = \calF(f(x))$.
\end{enumerate}
Given any sufficiently small $\epsilon>0$, there exists an $\epsilon$-tempered function $C(x)$ such that 
\begin{enumerate}[resume]
	\item there is  a measurable family of $C^r$-embeddings $\{\phi_x\colon \R^d(1)\to M\}$ with the following properties:
\begin{enumerate}
	\item $\phi_x(0) = x$;
	\item $\phi_x(\R^d(1))$ is a precompact open (in the immersed topology) neighborhood of $x$ in $\calF(x)$;
	\item $f(\phi_x(\R^d(1)))\subset \phi_{f(x)}(\R^d(1))$;
	\item $\frac 1 {C(x)}\le \|D\phi_x\|\le 1$.

\end{enumerate}
	\item \label{maps} If $\td f_x \colon \R^d(1) \to \R^d(1)$ is the map $$\td f_x = \phi_{f(x)} \inv \circ f \circ \phi_x$$
	then 
	\begin{enumerate}
		\item $\|\td f_x \|_{C^r}\le C(x)$
	\item  \label{contract} for all $v\in \R^d(1),$ $$\limsup _{n\to \infty} \frac 1 n \log \|\td f_{f^{n-1}(x)} \circ \dots \circ \td f_x(v)\| <0$$
\end{enumerate}
\item $\calF(x) = \bigcup _{n\ge 0}f^{-n} \left(\phi_{f^n(x)}(\R^d(1)) \right)$
\end{enumerate}
\end{definition}
Write  $$E_x = T_x \Fol(x)$$ and let $A^{(n)}(x) =\restrict{D_x f^n}{E_x}$. 
By \eqref{eq:integrability}, $x\mapsto \log \|A^{(1)}_x\|$ is $L^1(\mu)$ and condition  \eqref{contract} implies the top Lyapunov exponent for the dynamics tangent to $\Fol$ is negative: $\liminf_{n\to \infty} \frac 1 n \int \log  \|A^{(n)}_x\| \ d \mu(x) <0.$

Let $\Omega^+\subset M$ be the set of points $x\in M$ for which the sequence $\{A^{(n)}(x)\}$ is forwards regular.  It follows for $\mu$-a.e.\ $x\in \Omega^+$ that $\Fol(x)\subset \Omega^+$.  In particular, for a.e.\ $x$ and any $y\in \Fol(x)$, one obtains a forwards Lyapunov flag $\calV_y \subset E_y = T_y\Fol(x)$.  
It is well known (see e.g.\ \cite{MR556581}), that the variation $y\mapsto \calV_y$ is $C^{r-1}$ along $\Fol(x)$ for almost every $x$.

\subsection{Subresonant structures on contracting foliations}
We now suppose that $\mu$ is $f$-ergodic.  Let $\calF$ be a contracting, $C^r$-tame, $f$-invariant, measurable foliation.
Recall that for almost every $x\in \Omega_+$ and every $y\in \calF$, the sequence $\restrict{D_yf^n}{E_y}$ is forwards regular.  
Let
$$-\infty <\lambda_1<\dots<\lambda_\ell<0$$ be the Lyapunov exponents of $\restrict{D_xf^n}{E_x}$.

For every $x\in \Omega_+$, let $\varpi_x$ denote the Lyapunov {\SG} for the sequence $\{\restrict{D_xf^n}{E_x}\}$.
For $x\in \Omega_+$, let $\calV_x= \calV_{\varpi_x}$
 denote the associated forward Lyapunov flag in $E_x$.  

We will now assume $r>\lambda_1/\lambda_\ell.$

\subsubsection{Equivariant subresonant structures on leaves of contracting foliations}
The following follows from the main result \cite{MR3642250} after appropriate translations.  
\begin{theorem}[{c.f.\ \cite[Theorems 2.3 and 2.5]{MR3642250}}]\label{thm:normalforms}
There is a full measure subset $\Omega_0\subset \Omega+$  with $\Fol(x)\subset \Omega_0$ for every $x\in \Omega_0$ and a collection of $C^r$ diffeomorphisms $$\calH _x= \{ h\colon E_x \to \Fol (x)\}$$ with the following properties:
\begin{enumerate}
	\item {\bf (Finite dimensionality)} The group  $\calG^{SR}(E_x)$ of subresonant  polynomials acts transitively on $\calH_x$ be precomposition.  

	\item  \label{coher} {\bf (Coherence along leaves)}  For every $y \in \Fol(x)$, $\td h \in \calH _y$, and $ h\in \calH _x$, the transition function
	$$\td h\inv \circ  h\colon E_x\to E_y$$
	is a subresonant polynomial in $\calP^{SR}(E_x, E_y)$. 
	\item  {\bf (Equivariance under dynamics)} For $h\in \calH_x$ and $\hat h\in \calH_{f(x)}$,  the map $\hat h \inv \circ f\circ h\colon E_x\to E_{f(x)} $ is a  subresonant polynomial in $\calP^{SR}(E_x, E_{f(x)})$.
\item  {\bf (Measurability)} The assignment $x\mapsto \calH _x$ is measurable in the sense that it has a measurable section $x\mapsto h_x\in \calH_x$ with $h_x(0) = x$  and $D_0h_x  = \Id$.
\item {\bf (Uniqueness)} The families $\calH_x$ are uniquely defined $\mod 0$ by the above properties.
\end{enumerate}
Moreover, for any $\epsilon>0$ there exist 
tempered functions $\rho,C\colon M\to [1,\infty)$ and a measurable section  $x\mapsto h_x\in \calH_x$ such that $h_x(0) = x$, $D_0h_x  = \Id$,
and  $$\|\restrict{h_x}{E_x(\rho(x)\inv)}\|_{C^r}\le C(x).$$
\end{theorem}

We note that conclusion  \eqref{coher} is not stated in \cite[Theorem 2.3]{MR3642250} (which is formulated for non-linear extensions).  However, it is stated in \cite[Theorem 2.5(iv)]{MR3642250} and follows directly from our \cref{prop:forcedSR}.

\subsection{Fast foliations}
Given $\kappa<0$ and $x\in \Omega_x$, let $$E^\kappa_x:= \bigl\{v\in E_x: \limsup_{n\to \infty } \frac 1 n \log \|D_xf^n v\|\le \kappa\bigr\}$$
and write $$\Fol^{\kappa}(x)\subset \Fol(x):= \bigl\{ y\in \fol(x): \limsup_{n\to \infty} \frac 1 n \log d(f^n(x), f^n(y) ) \le \kappa \bigr\}$$
for the $\kappa$-fast submanifold of $\Fol(x)$.

\begin{claim}
For   any $\kappa<0$ and almost every $x\in \Omega_+$ and every $y\in \Fol(x)$
the set  $\Fol^{\kappa}(y)$ is an embedded $C^r$-submanifold  of $\Fol(x)$ tangent to $E^{\kappa}_x$.

If $r\ge 2$ the  collection of fast manifolds $\{\Fol^{\kappa}(y): y\in\Fol(x)\}$ forms a $C^r$-foliation of $\Fol(x)$. 
If $1<r<2$ and if $r>\lambda_1/\lambda_\ell$ then collection of fast manifolds $\{\Fol^{\kappa}(y): y\in\Fol(x)\}$ forms a $C^r$-foliation of $\Fol(x)$.
\end{claim}

Indeed, this follows from the existence of the normal form coordinate change in \cref{thm:normalforms}.

\subsubsection{Induced subresonant structures on $\kappa$-fast foliations}
Given $x\in \Omega_0$ and $\kappa<0$, recall we write $\calE^{\kappa}_x\subset \calE_x$,  $$\calE^{\kappa}_x= \{ v\in \calE_x: \varpi_x(v)\le \kappa\}.$$
We have the following which characterizes the fast stable manifolds $W^{\kappa}_y$ inside $W^\calE_x$ relative to subresonant structure on $W^\calE_x$ and shows that subresonant structures restrict on $W^{\calE}_x$ restrict to subresonant structures on $W^{\kappa}_x$.  
\begin{proposition}\label{prop:fast}
Given $x\in \Omega_0$  the following hold: 
\begin{enumerate}
	\item \label{fast1} Given $y\in\Fol(x)$, $v\in E_x$, and  $h \in \calH_x$ with $y= h (v)$ the restriction $$\restrict{h}{v+ E_x^{\kappa} }\colon v+ E_x^{\kappa}\to \Fol(x)$$ is a $C^r$ diffeomorphism between the affine subspace $v+ E_x^{\kappa}$ and $\Fol^{\kappa}(y)$. 
	\item The set of restrictions $$\calH^{\kappa}_x:=\{\restrict{h}{E_x^{\kappa} }:h\in \calH_x, h(0) \in  \Fol^{\kappa}_x\}$$
	defines a subresonant atlas on each $\Fol^\kappa_x$ and the collection $\{\calH^{\kappa}_x: x\in \Omega_0\}$ defines the unique equivariant subresonant structure on $\Fol^\kappa$.  
\end{enumerate}
\end{proposition}
\begin{corollary}\label{cor:restrictSR}

Fix $x\in \Omega_0$, $y\in\Fol(x)$,   $h_x \in \calH_x$, $h_y \in \calH_y$, and $v\in E_x$ with $y= h (v)$.  
Then the map $E_x^{\kappa} \to E_y^{\kappa}$,
$$u\mapsto h_y\inv (h_x(v+u))$$
is a subresonant polynomial.  
\end{corollary}
\subsection{Lyapunov flags and one-sided Lyapunov metric}
For $x\in M$, write $A^{(n)}(x) =\restrict{D_x f^n}{E_x}$.  By Oseldec's theorem, for almost every $x\in M$, the sequence  $\{A^{(n)}(x)\}$ is forwards regular.  Let $\lambda_1<\dots<\lambda_\ell$ be the Lyapunov exponents with associated multiplicities $m_1, \dots, m_\ell$. 

Using the charts $\phi_x$ as in \cref{def:contfol},  for almost every $x\in M$ we identify $E_x$ with $\R^d$ and equip $E_x$ with the  Euclidean norm $\|\cdot\|$ induced by this identification.  
 For $x\in M$, write $\varpi\colon E_x\to \R$ for the Lyapunov weight.

We have the following whose proof is nearly the same as \cref{absLyapu}.  
\begin{proposition}\label{prop:LN}
Fix $\epsilon>0$.  There  exists a measurable function $L_\epsilon\colon M\to (0,\infty)$ and for almost every  $x\in M$ an inner product $\iprod{\cdot,\cdot}'_x$ with associated norm $\|\cdot\|'_x$ on $E_x$ with the following properties: 
\begin{enumerate}
\item \label{cock1}$\|A^{(n)}_xv\|_{f^n(x)}'\le e^{n (\varpi(v)+\epsilon)}\|v\|_x'$  for every $n\ge 0$ and for every $v\in E_x$;
\item\label{cock2} $  \|v\| \le  \|v\|'_{f^n(x )}\le L_\epsilon (x)  e^{n (\frac{\epsilon}2)} \|v\|$ for all $v\in E_{f^n(x)}$;
\item  \label{cock3}$\|A^{(n)}_xv\|_{f^n(x)}'\ge e^{n (\varpi(v)-\epsilon)}\|v\|_x'$  for all $n$ sufficiently large (depending on $v$).
\end{enumerate}
\end{proposition}

\subsection{Centralizers}

Let $g\colon M\to M$ be a diffeomorphism that preserves the measure $\mu$.  Suppose  that $g$ commutes with $f$ and for almost every $x$,  $g(\calF(x)) = \calF(g(x))$.  
\begin{proposition}[{c.f.\ \cite[Theorem 2.3(3)]{MR3642250}}]
Suppose $g$ is $C^{r}$ for some $r>\lambda_1/\lambda_\ell$.  Then 
for $h\in \calH_x$ and $\hat h\in \calH_{g(x)}$,  the map $$\hat h \inv \circ g\circ h\colon E_x\to E_{g(x)}$$ is a  subresonant polynomial in $\calP^{SR}(E_x, E_{g(x)})$.
\end{proposition}
\begin{proof}
Pick a measurable section $x\mapsto h_x\colon E_x\to \Fol(x)$ with $h_x(0) = x$ in $\calH_x$.  For each $x$, set $$\phi_x= h_{g(x)}\inv\circ  g\circ h_{x}\colon E_x \to E_{g(x)}.$$
Fix $\epsilon>0$ sufficiently small and equip each $E_x$ with  the norm $\|\cdot \|'_x$ from \cref{prop:LN}.  
For $x\in M$ and $j \ge 0$, set  
\begin{enumerate}
\item $W^j_x =E_{f^j(x)} $
\item $V^j_x =  E_{f^j(g(x))}=  D_{f^j(x)} gE_{f^j(x)} $
\item $\td g^j_x \colon W^j_x\to W^{j+1}_x$, $\td g^j_x = h_{f^{j+1}(x)}\inv \circ  f\circ h_{f^j(x)}$
\item $\td f^j_x \colon V^j_x\to V^{j+1}_x$, $\td f^j_x = h_{f^{j+1}(g(x))}\inv \circ  f\circ h_{f^j(g(x))}$
\item $\td \phi_x^j\colon W^j_x\to V^j_x$,  $ \td \phi_x^j := \phi_{f^j(x)}$.
\end{enumerate}
Then, relative to the norms $\| \cdot \|'$ from \cref{prop:LN}, the maps $\td g^j_x $ and $\td f^j_x $ satisfy the hypotheses of  \cref{prop:forcedSR}.  Moreover, by Poincar\'e recurrence, for almost every $x$, relative to the norms  $\|\cdot \|'$ from \cref{prop:LN} 
the $C^r$ norm of the restriction of  $\td \phi^j_x $ to the unit ball is uniformly bounded for an infinite set of $j\in \N$.  
Moreover, $$\td \phi_x^n\circ \td g^{(n)}_x = \td f^{(n)}_x \circ \td \phi_x^0.$$
It then follows from  \cref{prop:forcedSR} that $\phi_x= \td \phi_x^0$ is a subresonant polynomial map.  
 \end{proof}

\subsection{Linearization of  dynamics along contracting foliations}
\begin{proposition}
For almost every $x\in \Omega_0$ there are  
\begin{enumerate}
\item	a vector space $\bfL \calF(x)$ and a codimension-1 subspace 
 $(\bfL \calF(x))'$;
\item an injective $C^r$ embedding $\iota_x\colon \calF(x) \to \bfL \calF(x)$ whose image is contained in the affine subspace $\iota_x(x) +  (\bfL \calF(x))'$;
\item a linear map $\bfL_x f\colon \bfL \calF(x) \to \bfL \calF(f(x))$ such that $$\bfL_x f\circ \iota_x = \iota_{f(x)} \circ f.$$
\end{enumerate}
Moreover the following hold:
\begin{enumerate}[resume]
\item For $y\in \calF(x)$ we have equalities	$\bfL \calF(y)= \bfL \calF(y)$, $\iota_x= \iota_y$, and $\bfL_x f= \bfL_y f$.
\item For every $h_x\in \calH_x$, there is an inner product on $\bfL \calF(x)$.  For any choice of inner product, the vector $\iota_x(x)$ is a unit vector orthogonal to  $(\bfL \calF(x))'$.
\item If $x\mapsto h_x$ is a tempered choice of coordinates then, relative to the induced inner products on $\bfL \calF(x)$, the Lyapunov exponents of the cocycle $\bfL _x f$ are all expressions of the form
$$ \lambda_1 \le  \sum_{i = 1} ^\ell n_i \lambda_i \le 0 $$
where $ n_i$ are non-negative integers. 

For the restriction to  $(\bfL \calF(x))'$, the Lyapunov exponents of the cocycle $\bfL _x f$ are all expressions of the form
$$ \lambda_1 \le  \sum_{i = 1} ^\ell n_i \lambda_i < 0 $$
where $ n_i$ are non-negative integers. 
\end{enumerate}

\end{proposition}


\begin{bibdiv}
\begin{biblist}

\bib{MR2348606}{book}{
      author={Barreira, Luis},
      author={Pesin, Yakov},
       title={Nonuniform hyperbolicity},
      series={Encyclopedia of Mathematics and its Applications},
   publisher={Cambridge University Press},
     address={Cambridge},
        date={2007},
      volume={115},
        ISBN={978-0-521-83258-8; 0-521-83258-6},
        note={Dynamics of systems with nonzero Lyapunov exponents},
      review={\MR{2348606 (2010c:37067)}},
}

\bib{MR2090767}{incollection}{
      author={Feres, R.},
       title={A differential-geometric view of normal forms of contractions},
        date={2004},
   booktitle={Modern dynamical systems and applications},
   publisher={Cambridge Univ. Press, Cambridge},
       pages={103\ndash 121},
      review={\MR{2090767}},
}

\bib{MR1908557}{article}{
      author={Guysinsky, M.},
       title={The theory of non-stationary normal forms},
        date={2002},
        ISSN={0143-3857,1469-4417},
     journal={Ergodic Theory Dynam. Systems},
      volume={22},
      number={3},
       pages={845\ndash 862},
         url={https://doi.org/10.1017/S0143385702000421},
      review={\MR{1908557}},
}

\bib{MR3642250}{article}{
      author={Kalinin, Boris},
      author={Sadovskaya, Victoria},
       title={Normal forms for non-uniform contractions},
        date={2017},
        ISSN={1930-5311},
     journal={J. Mod. Dyn.},
      volume={11},
       pages={341\ndash 368},
  url={https://doi-org.turing.library.northwestern.edu/10.3934/jmd.2017014},
      review={\MR{3642250}},
}

\bib{MR3893265}{article}{
      author={Melnick, Karin},
       title={Non-stationary smooth geometric structures for contracting
  measurable cocycles},
        date={2019},
        ISSN={0143-3857,1469-4417},
     journal={Ergodic Theory Dynam. Systems},
      volume={39},
      number={2},
       pages={392\ndash 424},
         url={https://doi.org/10.1017/etds.2017.38},
      review={\MR{3893265}},
}

\bib{MR556581}{article}{
      author={Ruelle, David},
       title={Ergodic theory of differentiable dynamical systems},
        date={1979},
        ISSN={0073-8301},
     journal={Inst. Hautes \'Etudes Sci. Publ. Math.},
      number={50},
       pages={27\ndash 58},
         url={http://www.numdam.org/item?id=PMIHES_1979__50__27_0},
      review={\MR{556581 (81f:58031)}},
}

\end{biblist}
\end{bibdiv}


\end{document}